\documentclass[11pt,twoside]{article}
\usepackage{amsfonts,amssymb,amscd,xy}
\usepackage[leqno]{amsmath}
\usepackage{mathrsfs}

\def\bm#1{\mathpalette\bmstyle{#1}}
\def\bmstyle#1#2{\mbox{\boldmath$#1#2$}}

\newcommand{\Cal}[1]{{\mathcal {#1}}}
\newcommand{\Lie}[1]{\mathrm {Lie}(#1)}
\newcommand{\uLie}[1]{{\underline{\mathrm{Lie}}}(#1)}

\newcommand{\id}{\mathrm{id}}

\newcommand{\dR}{\mathrm{dR}}

\newcommand{\mgr}[1]{{\mathbb G}_{m, #1}}
\newcommand{\agr}[1]{{\mathbb G}_{a, #1}}

\newcommand{\CC}{{\mathbb C}}
\newcommand{\ZZ}{{\mathbb Z}}

\newcommand{\TT}{{\mathrm T}}
\newcommand{\VV}{{\mathbb V}}

\newcommand{\spec}[1]{\mathrm {Spec}\left(#1 \right)}

\DeclareMathOperator{\coker}{coker}

\newcommand{\homo}[2]{\mathrm {Hom}(#1,#2)}
\newcommand{\uhomo}[2]{\underline{\mathrm {Hom}}(#1,#2)}
\newcommand{\uhomon}[2]{\underline{\mathrm {Hom}}^\nabla(#1,#2)}

\newcommand{\ext}[2]{{\mathrm {Ext }}^{1}(#1,#2)}

\newcommand{\extnat}[2]{{\mathrm {Ext }}^\natural(#1,#2)}
\newcommand{\biext}[3]{{\mathrm {Biext}}^1(#1,#2;#3)}
\newcommand{\uext}[2]{\underline{\mathrm{Ext}}^{1}(#1,#2)}
\newcommand{\uextnat}[2]{\underline{\mathrm {Ext }}^\natural(#1,#2)}

\newcommand{\uextfl}[2]{\underline{\mathrm {Ext}}^1_{\mathrm {fl}}(#1,#2)}
\newcommand{\uextzar}[2]{\underline{\mathrm{Ext}}^1_{\mathrm {Zar}}(#1,#2)}

\newcommand{\proof}{{\it Proof.\ }}
\def\squareforqed{\hbox{\rlap{$\sqcap$}$\sqcup$}}
\def\qed{\ifmmode\else\unskip\quad\fi\squareforqed}

\newcommand{\nobd}{\nobreakdash}

\newtheorem{lemma}{Lemma}[section]{\bf}{\it}
\newtheorem{theorem}[lemma]{Theorem}{\bf}{\it}
\newtheorem{proposition}[lemma]{Proposition}{\bf}{\it}
\newtheorem{corollary}[lemma]{Corollary}{\bf}{\it}
\newtheorem{definition}[lemma]{Definition}{\bf}{\rm}
\newtheorem{remark}[lemma]{Remark}{\bf}{\rm}
\newtheorem{example}[lemma]{Example}{\bf}{\rm}

\begin{document}

\input xy
\xyoption{all}

\title{Deligne's duality for de Rham realizations of 1-motives}
\author{Alessandra Bertapelle}
\maketitle

\begin{abstract}
We show that the pairing on de Rham realizations of 1-motives in "Theorie di
Hodge~III", IHES 44, can be defined over any base scheme and we prove that it
gives rise to a perfect duality if one is working with a 1-motive and its
Cartier dual. Furthermore, we study universal extensions of 
$1$\nobd-motives and their relation with $\natural$-extensions.
\end{abstract}


\section{Introduction}

It is known (cf. \cite{MM}) that the Lie algebra of the universal extension
$A^{\natural}$ of an abelian variety $A$ is canonically isomorphic to the
first de Rham cohomology group of the dual abelian variety $A'$.  Deligne
defines in \cite{De} the de Rham realization $\TT_\dR(M)$ of a $1$-motive
$M=[X\to G]$ over a base scheme $S$ as the Lie algebra of ${\bf G}^\natural$,
where $M^\natural=[X\to {\mathbf{G}}^\natural]$ is a universal extension of
$M$. In this way he gets a (covariant) functor from the category of $1$-motives to the
category of locally free sheaves over $S$. In the case $S$ is the spectrum of
an algebraically closed field $k$, Deligne defines a pairing 
$\Phi: \TT_\dR(M)\otimes \TT_\dR(M')\to k$ 
(cf. \cite{De}, 10.2.7.3) between the de Rham realizations of a $1$-motive and its
Cartier dual $M'=[X'\to G']$; for $k=\CC$ the pairing $\Phi$ coincides with an
analogous (perfect) pairing on Hodge realizations (cf. \cite{De}, 10.2.8) and
hence it is perfect. 
In the present paper we construct a pairing $\Phi$ between de Rham
realizations of dual $1$-motives over a general base $S$ and we show that it
is perfect; this fact generalizes also a result of Coleman that shows the
perfectness of $\Phi$ in the case of abelian schemes (over a base flat over
$\ZZ$) via the comparison with  a second (perfect) pairing (cf. \cite{Co},
1.1.1). 
As the existing proofs do not extend to the general case, we show 
directly the perfectness of $\Phi$ proving that this pairing  fits in a diagram
\begin{eqnarray*}
\xymatrix @C-=1pt{
{\boldsymbol{\omega}}_{G'}\ar[d]^{\iota}&\otimes&\uLie{G'}\ar[rr]&\quad &\uLie{\mgr S}\\
\Phi\colon \TT_\dR(M) \ar[d]^{g}& \otimes &
  \TT_\dR(M') \ar[u]^{g'}\ar[rr]& &\uLie{\mgr S}\\
\uLie{G} &\otimes&{\bm \omega}_{G}\ar[rr]\ar[u]^{\iota'}& &\uLie{\mgr S}}
\end{eqnarray*}
where the upper (resp.\ lower) pairing is the usual duality between
the Lie algebra of $G'$ (resp.\  $G$) and the sheaf of invariant
differentials of $G'$ (resp.\  $G$). As the maps $\iota$, $g'$  and
$\iota'$, $g$ come out to be transposes of each others, we get the perfectness
of $\Phi$ with no restriction on the base.

In the last section we describe the relation between a universal extension
$v\colon X\to {\mathbf{G}}^\natural$ of $M$ and  $\natural$\nobd-extensions of 
$M'$ showing that there is an exact sequence
\begin{eqnarray*}
\xymatrix{X\ar[r]^(.5){v}&{\mathbf{G}}^\natural\ar[r]&
\uextnat{M'}{\mgr S} \ar[r]&0. }
\end{eqnarray*}
This result generalizes the fact that the universal extension of an abelian
scheme $A$ represents the functor which assigns to an $S$-scheme $S'$ the 
$\natural$-extensions of the dual abelian scheme $A'$ by the multiplicative 
group over $S'$.
 

\section{Universal extensions of  $1$-motives}\label{sec.ue}

Let $S$ be a scheme. Recall that an \emph{$S$\nobd-$1$\nobd-motive}
$M=[u\colon X\to G]$ is a two term complex $($in degree $-1, 0)$
of commutative group schemes over $S$ such that $X$ is an $S$\nobd-group
scheme that locally for the \'etale topology on $S$ is isomorphic to a
constant group of type $\ZZ^r$, $G$ is an $S$\nobd-group scheme extension of an abelian
scheme $A$ over $S$ by a torus $T$,  $u$ is an $S$\nobd-homomorphism $X\to G$.
Morphisms of $S$\nobd-$1$\nobd-motives are usual morphisms of complexes.
The category of $1$-motives can be seen as a full subcategory of the
derived category of bounded complexes of fppf sheaves on $S$ (cf. \cite{RA}).

An \emph{extension} of an $S$\nobd-$1$\nobd-motive $M=[u\colon X\to
G]$ by a group $H$ is an extension $E$ of $G$ by $H$ together with a
homomorphism $v\colon X\to E$ that lifts $u$. Two extensions
$(E_i,v_i)$, $i=1,2$, are isomorphic if there exists an isomorphism
$\varphi\colon E_1\to E_2$ (as extension of $G$ by $H$) such that
$v_2=\varphi\circ v_1$. As usual, $\ext{M}{H}$ denotes the group of
isomorphism classes of extensions of $M$ by $H$. In the following,
we will simply speak of $1$\nobd-motives meaning $S$\nobd-$1$\nobd-motives.

A \emph{universal extension} of $M$ is an extension
$M^\natural=[X\to {\mathbf{G}}^\natural]$  of $M$
by a vector group $\VV(M)$ over $S$ such that the homomorphism of
push-out
\begin{eqnarray}\label{def.epsilon}
\epsilon\colon {\mathrm{Hom}}_{{\cal O}_S}(\VV(M),W)\longrightarrow \ext{M}{W}
\end{eqnarray}
 is an isomorphism for all vector groups
$W$ over $S$ (cf. \cite{De}). Observe that $M^\natural$ and $\VV(M)$
are determined up to canonical isomorphisms (cf. \cite{MM}, p. 2).
Universal extensions of $1$\nobd-motives
 exist (see \cite{De}, \cite{BVS}).
As  explained in \cite{MM}, I, 1.7, it is sufficient to show that
the following conditions are satisfied:
\begin{itemize}
\item[a)] $\uhomo{M}{\agr S}=0$, \item[b)] $\uextzar{M}{\agr S}$
is a locally free sheaf  of ${\cal O}_S$\nobd-modules of finite rank,
\end{itemize}
as sheaves for the Zariski topology over $S$. If this is the case,
\begin{eqnarray*}\uextzar{M}{{\cal W}}=\uextzar{M}{\agr S}
\otimes_{{\cal O}_S} {{\cal W}}
\end{eqnarray*} 
for any locally free ${\cal O}_S$\nobd-module of finite rank $\cal W$  
and one takes as $\VV(M)$ the vector group associated to the dual sheaf of 
$\uextzar{M}{\agr S}$.

For the torus $T$ and the abelian scheme $A$ condition a) is
automatically satisfied. Hence the same holds for the semi-abelian
scheme $G$ and then for $M$. As $\ext{T}{\agr S}=0$ also condition
b) holds for tori. For abelian varieties the result is proved in
\cite{MM}, 1.10. Moreover, denote by $A^\natural$ a (fixed) universal
extension of $A$; as $\ext{G}{\agr S}=\ext{A}{\agr S}$ one gets
$\VV(G)=\VV(A)$ and a universal extension of $G$ is
$G^\natural=A^\natural\times_A G$ (see \cite{ABV}, 2.2.1). 
\begin{lemma}\label{lem.sheaf}
Let $M=[u\colon X\to G]$ be a $1$\nobd-motive and $W$  a vector
group over $S$. Then the functor
\begin{eqnarray}\label{funct.MW}
S'\leadsto \ext{M_{S'}}{W_{S'}}
\end{eqnarray}
is a sheaf for the flat and Zariski topologies. 
Here $M_{S'}$ denotes the $S'$\nobd-$1$\nobd-motive obtained via base-change.
\end{lemma}
\proof Consider the sequence
\begin{eqnarray}\label{seq.motive}
\xymatrix{0\ar[r]& G\ar[r]&M\ar[r] & [X\to 0]\ar[r]&0 }
\end{eqnarray} 
and recall that the functor  $S'\leadsto
\homo{X_{S'}}{Q_{S'}}$ is a sheaf for any $S$\nobd-group scheme $Q$
and that the category $\mathrm{EXT}(G,W)$ is rigid
(cf.~\cite{MM}, I, 1.10 proof). \qed

From (\ref{seq.motive}) we get an exact sequence
\begin{eqnarray}\label{seq.vectorial1}
\xymatrix{0\ar[r]&X^*=\uhomo{X}{\agr S}\ar[r]&\uextzar{M}{\agr
S}\ar[r] & \uextzar{G}{\agr S} \ar[r] &0 }
\end{eqnarray}
that might not be exact on the right. However, it is exact on a suitable 
affine \'etale covering of $S$ where $X$ becomes constant. Since
 $X^*$ and $\uextzar{G}{\agr S}$ are locally free of finite
rank, the same is $\uextzar{M}{\agr S}$ and (\ref{seq.vectorial1}) is exact 
on the right. In particular,
$M$ admits a universal extension $M^\natural$. Passing to duals on
(\ref{seq.vectorial1}) one gets also a sequence of vector groups
\begin{eqnarray}\label{seq.vectorial2}
 \xymatrix {0\ar[r]&\VV(G)\ar[r]^(.5)i & \VV(M) \ar[r]^(.4){\bar \tau} & X\otimes \agr
S\ar[r] &0. } 
\end{eqnarray}

\subsection{A description of $\VV(M)$ via invariant differentials.}
It is well known that given an abelian scheme $A$ over $S$ the
vector group  $\VV(A)$ corresponds to the locally free sheaf
$\uLie{A'}^*={\boldsymbol{\omega}}_{A'}$ of invariant differentials of the
dual abelian scheme $A'$. In the next pages, we will generalize this
result to $1$\nobd-motives showing that if  $M'=[u'\colon X'\to G']$
is the Cartier dual of $M$, the vector group
 $\VV(M)$ corresponds to the sheaf ${\boldsymbol{\omega}}_{G'}$ of invariant 
differentials of the semi-abelian scheme $G'$. This fact will be of great use 
in the following sections.

For the definition of the dual motive $M'=[u'\colon X'\to G']$ of
$M=[u\colon X\to G]$ we refer to \cite{De}. Denote by
$[X\to A]$ the $1$\nobd-motive obtained via composition of $u$ with
the homomorphism $G\to A$. We recall that by definition  $G'$
represents the sheaf $\uextfl{[X\to A]}{\mgr S}$, the group $X'$ is
the group  of characters of $T$  and $u'$ is the boundary
homomorphism of the long exact sequence of Ext sheaves obtained
applying $\underline{\mathrm{Hom}}(-, \mgr S)$ to the exact
sequence
\begin{eqnarray}\label{seq.motiveMMA}
  \xymatrix{ 0\ar[r]& T\ar[r]&M\ar[r] & [X\to A] \ar[r] &0 }.
\end{eqnarray}
Furthermore, the sequence
\begin{eqnarray}\label{seq.motiveMA}
\xymatrix{
 0\ar[r]& A \ar[r]& [X\to A]\ar[r]& [X\to 0] \ar[r]& 0
 }
\end{eqnarray}
provides a short exact sequence
\begin{eqnarray}\label{seq.Gdual}
\xymatrix @C=13pt{
 0\ar[r]&T'=\uhomo{X}{\mgr S}\ar[r]&G'=
 \uextfl{[X\to A]}{\mgr S}\ar[r] &A'=
 \uextfl{A}{\mgr S} \ar[r] &0 }
 \end{eqnarray}
that describes $G'$ as a semi-abelian scheme. In the case $A=0$ and
$M=[u\colon X\to T]$ the dual $1$\nobd-motive $M'=[u'\colon X'\to
T']$ is simply obtained via the usual Cartier duality.

We start relating the Lie algebra of a semi-abelian scheme to vector
extensions of its Cartier dual.

\begin{lemma}\label{lem.liesemiab}
Let  $B$ be a semi-abelian scheme over $S$. Then
\begin{eqnarray*}\uLie{B}=\uLie{\uextfl{N}{\mgr S}}=
\uextfl{N}{\agr S}\end{eqnarray*}
where $N$ is the  $1$\nobd-motive Cartier dual of $[0\to B]$. In
particular, if we think $\uLie{B}$  as a sheaf for the Zariski
topology then
\begin{eqnarray*}\uLie{B}=\uextzar{N}{\agr S}.\end{eqnarray*}
\end{lemma}
\proof The last assertion follows from the first via Lemma~\ref{lem.sheaf}. 
The first isomorphism is obvious because $B$ is
isomorphic to $\uextfl{N}{\mgr S}$ by Cartier duality. It remains to
prove that $\uLie{B}$ is isomorphic to $\uextfl{N}{\agr S}$.

Given a scheme $S'$, denote by $S_\epsilon'$ the fibre product
$S'\times_{\spec \ZZ} \spec{\ZZ[\epsilon]/(\epsilon^2) }$. Recall that
by definition of Lie algebras we have an exact sequence
\begin{eqnarray}\label{eq.lie}
\xymatrix
{\quad\quad
 0\ar[r]& \Lie{B_{S'}} \ar[r]&
B(S'_\epsilon)=\ext {N_{S'_\epsilon}}{\mgr {S_\epsilon'}}
\ar[r]^(.5){f_B}& B(S')=\ext{N_{S'}}{\mgr {S'}}
 }
\end{eqnarray}
where $f_B$ is the composition with the closed immersion $S'\to
S_\epsilon'$ induced by
 \begin{eqnarray*}\ZZ[\epsilon]/(\epsilon^2) \to \ZZ,
 \quad \epsilon \mapsto 0,
 \end{eqnarray*}
(or the base-change on exact sequences). Let now
$\Re_{S_\epsilon'/S'}(\mgr {S_\epsilon'} )$ be the Weil restriction
of $\mgr {S_\epsilon'}$ with respect to the base-change morphism
$S'_\epsilon\to S'$. For any $S'$\nobd-scheme $Z$ it holds
$\Re_{S_\epsilon'/S'}(\mgr {S_\epsilon'} )(Z)= \mgr
{S_\epsilon'}(Z_\epsilon)=\mgr {S'}(Z_\epsilon)$. Moreover, there is
an exact sequence
\begin{eqnarray}\label{seq.weil}
\xymatrix{
 0\ar[r]& \agr {S'} \ar[r]&  \Re_{S_\epsilon'/S'}(\mgr {S_\epsilon'} )
 \ar[r]^(.6){f}& \mgr {S'} \ar[r]& 0
 }
\end{eqnarray}
 where the homomorphism $f\colon \mgr {S'} (Z_\epsilon)
\to \mgr {S'}(Z)$  is obtained via composition with the closed
immersion $Z\to Z_\epsilon$. From (\ref{seq.weil}) we deduce an
exact sequence
\begin{eqnarray}\label{eq.kerfN}
\xymatrix{
 0\ar[r]&  \ext{N_{S'}}{\agr {S'}} \ar[r]&
\ext{N_{S'}}{\Re_{S_\epsilon'/S'}(\mgr {S_\epsilon'} ) }
\ar[r]^(.6){f_N}& \ext{N_{S'}}{\mgr {S'}}
 }
\end{eqnarray}
where $f_N$ is the push-out with respect to $f$.

 In order to prove that $\uLie{B}=\uextfl{N}{\agr S}$, it is
sufficient to check  that
\begin{eqnarray*}
\Lie{B_{S'}}=\ext {N_{S'}}{\agr {S'}}
\end{eqnarray*}
for any $S$-scheme  $S'$. Comparing (\ref{eq.lie}) and
(\ref{eq.kerfN})
 we are reduced to see that
\begin{eqnarray*}
\ext{N_{S'}}{\Re_{S_\epsilon'/S'}(\mgr {S_\epsilon'} ) }=
\ext{N_{S'_{\epsilon}}}{\mgr {S_\epsilon'}  }
\end{eqnarray*}
and that $f_N$ coincides with $f_B$. This is not hard using
properties of Weil restriction. \qed

Recall that $\VV(M)$ is the vector group associated to
the dual  sheaf of $\uextzar{M}{\agr S}$. We will prove now that it
corresponds to the sheaf of invariant differentials of $G'$.

\begin{proposition} \label{pro.VMlie}
Let $M$ be a $1$\nobd-motive. It holds
\begin{eqnarray}
    \uext{M}{\agr S}= \uext{[X\to A]}{\agr S}= \uLie{G'}
\end{eqnarray}
for the flat and Zariski topologies. Hence $\VV(M)$ is (the
vector group associated to) the sheaf of invariant differentials
${\boldsymbol{\omega}}_{ G'}$.
 Moreover, the sequence of
vector groups in $(\ref{seq.vectorial2})$ is  the sequence
\begin{eqnarray}\label{seq.invdiff}
\xymatrix{
    0\ar[r]&{\boldsymbol{\omega}}_{A'} \ar[r]&{\boldsymbol{\omega}}_{G'}  \ar[r] &
    {\boldsymbol{\omega}}_{T'} \ar[r] &0 }
\end{eqnarray}
of invariant differentials of $(\ref{seq.Gdual})$.
\end{proposition}

\proof Denote by $M_A$ the $1$\nobd-motive $[X\to A]$. The first
isomorphism $\uext{M}{\agr S}= \uext{M_A}{\agr S}$ comes from the
 exact sequence in (\ref{seq.motiveMMA}) using the vanishing
$\uhomo{T}{\agr S}= 0=\uext{T}{\agr S}$. The second isomorphism was
proved in the previous lemma for $B=G'$ and $N=M_A$.

 For the second assertion, observe that using the isomorphisms
\begin{eqnarray*}\uextzar{M}{\agr S}= \uextzar{M_A}{\agr S}, \quad
\uextzar{G}{\agr S}= \uextzar{A}{\agr S},
\end{eqnarray*}
 the sequence (\ref{seq.vectorial1}) coincides with the sequence
\begin{eqnarray}
\label{seq.vectzar}
\xymatrix{  0\ar[r]&X^*=\uhomo{X}{\agr S}\ar[r]&\uextzar{M_A}{\agr
S}\ar[r] & \uextzar{A}{\agr S} \ar[r] &0 }
\end{eqnarray}
obtained from (\ref{seq.motiveMA}).
Now, the proof of Lemma~\ref{lem.liesemiab}  says that
(\ref{seq.vectzar}) is the sequence of Lie algebras
\begin{eqnarray*}
\xymatrix{
    0\ar[r]&\uLie{T'} \ar[r]&\uLie {G'}  \ar[r] &
    \uLie {A'} \ar[r] &0  }
\end{eqnarray*}
of the sequence (\ref{seq.Gdual}). Passing to duals we get the
desired result. \qed

\subsection{A description of ${\bf G}^\natural$ as push-out}

By the universal property of the universal extension of $G$ and 
Proposition~\ref{pro.VMlie}, the group scheme 
${\mathbf{G}}^\natural$ is (isomorphic to) the push-out
\begin{eqnarray}\label{dia.vectorial}
\xymatrix{
 0\ar[r]&\VV(G)={\boldsymbol{\omega}}_{A'} \ar[r]\ar[d]^i&G^\natural
     \ar[r]^(.5){\rho}\ar[d] & G\ar[r]\ar@2{-}[d] &0\\
0\ar[r]&\VV(M)={\boldsymbol{\omega}}_{G'}\ar[r]^(.6)\iota
\ar[d]^{\bar\tau}&{\mathbf{G}}^\natural \ar[r]\ar[d]^\tau & G\ar[r]
&0\\
& X\otimes \agr S={\boldsymbol{\omega}}_{T'}\ar@2{-}[r]&{\boldsymbol{\omega}}_{T'}&
&}\end{eqnarray} 
where the vertical sequence on the left is
(\ref{seq.vectorial2}) or (\ref{seq.invdiff}); this fact was firstly observed
in \cite{BVS} (without the contribution of invariant differentials).
There is then a useful criterion to test when a
homomorphism $X\to {\mathbf{G}}^\natural$ provides a universal extension of $M$:
\begin{lemma}\label{lem.criterion}
Let $M=[u\colon X\to G]$ be a $1$\nobd-motive as above and let
${\mathbf{G}}^\natural$ be the group scheme defined in diagram
(\ref{dia.vectorial}). A homomorphism $v\colon X\to
{\mathbf{G}}^\natural$ such that $\rho\circ v=u$ is a universal
extension of $M$ if and only if $\tau\circ v\colon X\to X\otimes
\agr S$ is a universal extension of the $1$\nobd-motive $[X\to 0]$.
\end{lemma}
\proof
 Let $W$ be a vector group over $S$ and consider the following
diagram
\begin{eqnarray*}\label{dia.criterion}
\xymatrix{
 0\ar[r]& {\mathrm{Hom}}_{{\cal O}_S}(X\otimes \agr
     S,W)\ar[r]^(.5){\bar\tau^*}\ar[d]^{\epsilon_X}&{\mathrm{Hom}}_{{\cal
           O}_S}(\VV(M),W)\ar[r]^(.5){i^*}\ar[d]^\epsilon & 
                 {\mathrm{Hom}}_{{\cal O}_S}(\VV(G),W) \ar[d]^\wr  \\
 0\ar[r]& \ext{[X\to 0]}{W}\ar[r] &\ext{M}{W} \ar[r]&
\ext{G}{W}. & 
}
\end{eqnarray*}
where the upper sequence is obtained from (\ref{seq.vectorial2}) and
the lower sequence is obtained from (\ref{seq.motive}).   Given a
morphism of vector groups  $f\colon X\otimes \agr S\to W$, 
 $\epsilon_X(f)$ is   the trivial extension of $0$ by  $W$ together
with the morphism $f\circ \tau \circ v\colon X\to W$; for a
$g\colon \VV(M)\to W$ the extension $\epsilon(g)$ is the
isomorphism class of the push-out with respect to $g$ of the
extension $v\colon X\to {\mathbf{G}}^\natural$ of $M$
 by $\VV(M)$.
The push-out homomorphism on the right is an isomorphism because of
the the universal property of $G^\natural$. By construction the diagram
is commutative.

If $v\colon X\to {\mathbf{G}}^\natural$ is a universal extension of $M$, the
homomorphism $\epsilon$  is an isomorphism and hence also
$\epsilon_X$ is an isomorphism. This says that $\tau \circ v$ is a
universal extension of $[X\to 0]$.
Suppose now that $\tau \circ v$ is a universal extension of $[X\to
0]$, i.e. $\epsilon_X$ is an isomorphism. This implies that
$\epsilon$ is injective. If ${i^*}$ is surjective, we can deduce
that also $\epsilon$ is an isomorphism and hence $v$ is a universal
extension of $M$. In general ${i^*}$ is surjective Zariski locally
on $S$. As the functor in (\ref{funct.MW})
is a Zariski-sheaf, as well as 
$U\leadsto {\mathrm{Hom}}_{{\cal O}_{S|U}}(\VV(M)_{|U},W_{|U})$, 
  we conclude that 
$\epsilon$ is  an isomorphism.
 \qed

\begin{remark}\label{rem.universal}
Observe that if $v\colon X\to {\mathbf{G}}^\natural $ is a universal 
extension of $M$ and $f\colon X\to \VV(M)$ is a homomorphism 
that factors through $\VV(G)$ then $v+f\colon X\to {\mathbf{G}}^\natural $ 
is a universal extension too. 
However, $v+f$ is isomorphic to $v$, as extension of $M$ by $\VV(M)$,  if and
only if $f=0$ because the extension ${\mathbf{G}}^\natural$ admits no 
non-trivial automorphisms.
\end{remark}

\section{$\natural$-structures}

In order to define Deligne's pairing for the de Rham
realizations of $1$\nobd-motives, we need  to recall first some
definitions and results on $\natural$\nobd-extensions and
$\natural$\nobd-biextensions. Proposition~\ref{pro.repE} of this section is
the key result that permits to generalize [5], 10.2.7.4.


\subsection{Some definitions}
Let $S$ be a fixed  scheme, $Z$ an $S$-scheme, $G,Y_1, Y_2,H_1,H_2$
commutative $S$\nobd-group schemes with $G$ smooth. As usual 
$G_Z$ means $G\times_S Z$. Denote by $\Delta^1(Z)$
the first infinitesimal neighborhood of the diagonal $Z\to Z\times_S
Z$ and by $p_j\colon \Delta^1(Z)\to Z$, $j=1,2$, the morphisms
induced by the  usual projections $p_j\colon Z\times_S Z \to Z$.

\begin{definition}[\cite{MM}] A $\natural$\nobd-$G$\nobd-\emph{torsor}
on $Z$ is a torsor (for the \'etale
topology\footnote{Cf. \cite{MM}, I, 3.1. This hypothesis is needed
 to defined the curvature form of a connection by descent.})
$P$ on $Z$ under $G_Z$ endowed with an integrable connection, {i.e.}
an isomorphism $\nabla\colon p_1^*P\to p_2^*P$ of
$G_{\Delta^1(Z)}$\nobd-torsors which restricts
to the identity on $Z$ and has zero curvature.
\end{definition}

The trivial  $\natural$\nobd-$G$-torsor is the trivial torsor $G_Z$
endowed with the trivial connection $\nabla^0$, {i.e.} the identity
on $G_{\Delta^1(Z)}$. A trivialization of a $\natural$\nobd-$G$-torsor
$(P,\nabla)$ is a section $s\colon Z\to P$ such that the induced
isomorphism $\varphi_s\colon (P,\nabla)\to (G_Z,\nabla^0)$ is
horizontal. Observe that given a trivialization $s$ of a torsor $P$
on $Z$ under $G$ there is a unique possible
$\natural$\nobd-structure that makes $\varphi_s$ horizontal. We
describe in details  a case that will be needed later.

\begin{example}\label{ex.1}
 Let $Z=\mgr S =G$, $P=G\times Z$ with the trivialization 
$s\colon Z\to P, ~b\mapsto (b^n,b)$. We have
an isomorphism $\varphi_s\colon P\to \mgr S^2$, $(x,y)\mapsto
(x/y^n, y)$ such that $\varphi_s\circ s$ is the usual trivialization
$b\mapsto (1,b)$. Consider a connection $\nabla$ on $P$ given by a
global differential $\omega$ on $Z$.  There is a unique possible
choice of  $\nabla$ that makes $\varphi_s$ horizontal with respect
to $\nabla$ on $P$ and the trivial connection $\nabla^0$ on $\mgr
S^2$. More precisely, let $t$ (resp.\  $z$) be the parameter of $G$
(resp.\  $Z$); the isomorphism $\varphi_s$ induces an isomorphism of
 algebras such that $\varphi^*_s(z)=z$ and $\varphi^*_s(t)=t/z^n$, hence an
isomorphism  $\varphi^*_s\colon  {\cal O}_Z\to {\cal
O}_Z, $ $ a\mapsto a/z^n $. The horizontality condition says that the
induced connections on sheaves (see \cite{MM}, I, 3.1.2) $\nabla,
\nabla^0\colon {\cal O}_Z\to \Omega_{Z/S}^1$ satisfy
\begin{eqnarray*}
0=\varphi^*_s(\nabla^0(1))=\nabla(\varphi_s^*(1))=
\nabla(1/z^n)=-ndz/z^{n+1}+(1/z^n)w
\end{eqnarray*}
and hence $w=ndz/z$.
\end{example}

Let in the following $Z$ be  a group scheme over $S$ and
denote by $\mu_Z\colon Z\times_S Z\to Z$ its group law.

\begin{definition}[\cite{MM}]\label{def.naturalext}
A $\natural$\nobd-\emph{extension} of $Z$ by $G$ is a
$\natural$\nobd-$G$\nobd-torsor $(P,\nabla)$ on $Z$ where $P$ is
 extension of $Z$ by $G$ and the usual morphism
\begin{eqnarray*}\nu\colon  p_1^*P+p_2^*P\to \mu_Z^* P
\end{eqnarray*}
is horizontal.
\end{definition}

Two $\natural$\nobd-extensions $(P_i,\nabla_i)$, $i=1,2$, of $Z$ by
$G$ are isomorphic if there exists an isomorphism (of extensions)
$\phi\colon P_1\to P_2$ that is horizontal. The trivial
$\natural$\nobd-extension of $Z$ by $G$ is the trivial extension
$Z\times_S G$ equipped with the trivial connection $\nabla^0$. A
trivialization of a $\natural$\nobd-extension $(P,\nabla)$ is a
section $s\colon Z\to P$ that provides an isomorphism of
$(P,\nabla)$ with the trivial $\natural$\nobd-extension. One denotes
by $\extnat{Z}{G}$ the group of isomorphism classes of
$\natural$\nobd-extensions of $Z$ by $G$. We have an exact sequence
(cf. \cite{MM}, II 4.2)
\begin{equation}\label{seq.homextnatext} 
\xymatrix{ 
\homo{Z}{G}\ar[r]&
\Gamma({\bm{\omega}_Z}\otimes \uLie{G})\ar[r]& \extnat{Z}{G}
\ar[r]^(.5)f & \ext{Z}{G}.}
\end{equation}

 Following  \cite{Co2}, 0.2, it is easy to describe
  $\natural$\nobd-extensions of $Z$ by $\mgr S$.

\begin{proposition}\label{pro.extnat}
Let $Z$ be a commutative group scheme over $S$ and $E$ an extension
of $Z$ by $\mgr S$. There is a one-to-one correspondence between
connections $\nabla$ on $E$ making $(E,\nabla)$  a
$\natural$\nobd-extension of $Z$ by $\mgr S$ and
\emph{normal}\footnote{Denote by $z$ the standard parameter of $\mgr
S$. An invariant differential on $E$ is said to be normal if it
pulls back to $dz/z$ on $\mgr S$.} invariant differentials on $E$.
\end{proposition}
\proof It is known (see \cite{Co2}, 0.2.1) that there is a
one-to-one correspondence between connections on $E$ and global
differentials on $E$ that pull back to $dz/z$ and are invariant
under the action of $\mgr S$. Now, the horizontality condition in
Definition~\ref{def.naturalext} requires the global differential of
$E$ to be  invariant. \qed

 We recall now some definitions from \cite{De}.

\begin{definition}\label{def.natbiext}
 Let $P$ be a biextension
 of $(H_1, H_2)$ by $G$ and consider the usual morphisms
\begin{eqnarray*}\label{def.nu}
\nu_1\colon  p_{13}^*P+p_{23}^*P\to (\mu_1\times \id)^* P ,
\quad \text{on} \quad  H_1\times_S H_1\times_S H_2 ,\\
\nu_2\colon  p_{12}^*P+p_{13}^*P\to (\id \times \mu_2)^* P
   \quad \text{on} \quad  H_1\times_S H_2\times_S H_2.
\end{eqnarray*}Here $p_{ij}$ are the obvious
projections and $\mu_i$ is the group law on $H_i$. A
$\natural$\nobd-\emph{structure} on $P$ is a connection $\nabla$ on
the $G$\nobd-torsor $P$  over $H_1\times H_2$ such that
$\nu_1,\nu_2$ are horizontal. We will also say that $(P,\nabla)$ is
a $\natural$-\emph{biextension} of $(H_1, H_2)$ by $G$.
\end{definition}

A trivialization of a $\natural$\nobd-biextension $(P,\nabla)$ is a
horizontal isomorphism of $P$ with the trivial biextension endowed
with the trivial connection.

\begin{definition}\label{def.12structure}
Let $P$ be a biextension of $(H_1, H_2)$ by $G$. A
$\natural$\nobd-\emph{1\nobd-structure} on $P$ is  a connection on
$P$ such that $P$ becomes a $\natural$\nobd-extension\footnote{Here
$H_2$ is seen as base scheme; $\nu_1$ is automatically horizontal
because of Definition~\ref{def.naturalext}.} of $H_{1,H_2}$ by
$G_{H_2}$ with $\nu_2$ horizontal. A
$\natural$-\nobd\emph{2\nobd-structure} on $P$ is  a connection on
$P$ such that $P$ becomes a $\natural$\nobd-extension
 of $H_{2,H_1}$ by $G_{H_1}$ with
$\nu_1$ horizontal.
\end{definition}
Giving  a $\natural$\nobd-structure to a biextension is equivalent
to giving a $\natural$\nobd-$1$\nobd-structure and a
$\natural$\nobd-$2$\nobd-structure.

\begin{definition}
A \emph{biextension of complexes} $([Y_1\to H_1]$, $[Y_2\to H_2])$
by $G$ is a biextension $P$ of $(H_1, H_2)$ by $G$ endowed with a
trivialization of the pull-back of $P$ to $Y_1\times_S H_2$ and a
trivialization of the pull-back of $P$ to $H_1\times_S Y_2$ that
coincide on $Y_1\times_S Y_2$.

A $\natural$\nobd-\emph{extension} of a complex $[u\colon Y\to H]$
by a group $G$ is a $\natural$\nobd-extension $(P,\nabla)$ of $H$ by
$G$ with a trivialization (as $\natural$\nobd-biextension) of
the pull-back of $(P,\nabla)$ to $Y$.

A $\natural$\nobd-\emph{biextension of complexes} $([Y_1\to H_1]$,
$[Y_2\to H_2])$  by $G$ is a $\natural$\nobd-biextension
$(P,\nabla)$ of $(H_1, H_2)$ by $G$ endowed with a trivialization
(as $\natural$\nobd-biextension)
 of the pull-back of $(P,\nabla)$ to
$Y_1\times_S H_2$  and a trivialization of the pull-back of
$(P,\nabla)$ to $ H_1\times_S Y_2$ that coincide on $Y_1\times_S
Y_2$.
\end{definition}

\subsection{$\natural$-structures and biextensions.}

It is shown in \cite{MM} that the universal extension $A^\natural$ of
an abelian scheme over $S$ represents the functor that associates to
any $S$\nobd-scheme  $S'$ the group  of isomorphism classes of
$\natural$\nobd-extensions of $A'_{S'}$ by $\mgr S'$. See also
\cite{Co2}, 0.3.1. We will prove in  Lemma~\ref{lem.extnatMA}
that $G^\natural=\uextnat{[X'\to A']}{\mgr S}$, or the same, that
$G^\natural$ represents  the pre-sheaf for the flat topology
\begin{equation}\label{functorEG}
  S' \leadsto \left\{
 \begin{array}{ll}
 (g,\nabla),& g\in G(S'), \nabla\text{ a }
\natural \text{-structure on the extension }\Cal P_g' \text{ of }\\
&[X'\to A'] \text{ by }\mgr {S'} \text{ associated to } g
\end{array}\right\}.
\end{equation}
Observe that ${\cal P}_g'$ is the fibre at $g$ of
  the Poincar\'e biextension  ${\cal P'}$ of $(G,[X'\to
A'])$ by $\mgr S$.
\smallskip

We can generalize the result above to any $1$\nobd-motive:

\begin{proposition}\label{pro.repE}
Let $M$ be a $1$\nobd-motive, $M'$ its Cartier dual and   $\cal P$
the Poincar\'e  biextension of $(M,M')$.
The  group scheme  ${\mathbf{G}}^\natural$  defined in (\ref{dia.vectorial})
represents  the pre-sheaf for the flat topology
\begin{equation*}
{\cal E}\colon S' \leadsto  \left\{
 \begin{array}{ll}
 (g,\nabla),& g\in G(S'), \nabla\text{ a }
\natural \text{-structure on the extension }\Cal P_g \text{ of }\\
&M' \text{ by } \mgr {S'} \text{ associated to } g
\end{array}
\right\}.
\end{equation*}
\end{proposition}
 Observe that $\Cal P_g$ is the fibre at $g$ of the Poincar\'e
biextension  $\cal P$ of $(M, M')$. The biextension  $\cal P$ is
also the pull-back of ${\cal P'}$ (the Poincar\'e biextension of
$(G,[X'\to A'])$) to $(G,M')$ together with a suitable
trivialization on $X\times G'$. Hence $\Cal P_g$ can be seen as the
pull-back to $M'$ of  $\Cal P_g'$. In the following we will denote
by ${\cal P}$ (resp.\  ${\cal P'}$) also the $\mgr S$\nobd-torsor over
$G\times G'$ (resp.\  over $G\times A'$) underlying ${\cal P}$
(resp.\  ${\cal P'}$). In particular, the fibre ${\cal P}_g$ at a
point $g\in G(S')$ can be read as the pull-back to $G'$ of the fibre
of ${\cal P'}$ at $g$:
\begin{eqnarray}\label{dia.Poinc}
\xymatrix {
 0\ar[r]&\mgr {S'}\ar[r]\ar@{=}[d] &{\cal P}_g \ar[r]\ar[d]
 &G'_{S'}\ar[d]
 \ar[r]  &0\\
 0\ar[r]&\mgr {S'}\ar[r] &{\cal P}_g' \ar[r]
 &A'_{S'}
 \ar[r]  &0.
 }
 \end{eqnarray}

\proof We will construct a canonical isomorphism
\begin{eqnarray}\label{def.psi}
\Psi\colon {\mathbf{G}}^\natural(S')\to {\cal E}(S').
\end{eqnarray}
It is immediate to show that $\cal E$ is a sheaf for the flat
topology. Hence, we reduce the proof to the case where $S'=S$ is
affine and the short exact sequence in (\ref{seq.invdiff})
\begin{eqnarray*}
\xymatrix {
 0\ar[r]&{\boldsymbol{\omega}}_{A'}\ar[r] &{\boldsymbol{\omega}}_{G'}\ar[r]^(.5){\bar \tau} &
{\boldsymbol{\omega}}_{T'}\ar@{.>}@/^/[l]^\delta   \ar[r]  &0
 }
 \end{eqnarray*}
is split over $S$. Recall the notations in (\ref{dia.vectorial}) and
that ${\mathbf{G}}^\natural$ is extension of $G$ by ${\boldsymbol{\omega}}_{G'}=\VV(M)$.
Let $\delta$ be a section of $\bar \tau$ and denote by $\delta$ also
the induced section of $\tau\colon {\mathbf{G}}^\natural\to {\bm
\omega}_{T'}=X\otimes \agr S$. We identify then ${\mathbf{G}}^\natural $ with
$G^\natural\oplus {\boldsymbol{\omega}}_{T'}$. An $S$\nobd-valued point $a$ of
${\mathbf{G}}^\natural$ becomes a sum $a-\delta(\tau(a))\oplus \tau(a)$ where
$a-\delta(\tau(a))\in G^\natural(S)$ corresponds to a
$\natural$\nobd-structure on ${\cal P}'_{\rho(a)}$ via the functor
in (\ref{functorEG}); let $\eta_{A,a}$ be the corresponding normal
invariant differential of ${\cal P}'_{\rho(a)}$ 
(cf. Proposition~\ref{pro.extnat}). 
Then $\eta_a:=\eta_{A,a}+\delta(\tau(a))$ is a
normal invariant differential of ${\cal P}_{\rho(a)}$ and hence it
provides a $\natural$\nobd-structure $\nabla_{\eta_a}$ on ${\cal
P}_{\rho(a)}$. Observe that the definition of $\eta_a$ makes sense
because of diagram (\ref{dia.Poinc}). We define then
$\Psi(a)=(\rho(a),\nabla_{\eta_a})$.

For the injectivity of $\Psi$, let $a,b$ be two $S$\nobd-valued
points of $ {\mathbf{G}}^\natural$ such that $({\rho(a)},\nabla_{\eta_a})= $
$({\rho(b)},\nabla_{\eta_b} )$. Then
 $\rho(a)=\rho(b)$ and $\eta_a=\eta_b$. Define
$\omega:=a-b\in {\boldsymbol{\omega}}_{G'}(S)=(\ker \rho)(S)$. It holds
\begin{eqnarray*}
\eta_a&=&\eta_{A,a}+\delta(\tau(a))=\eta_{A,b}+\omega-\delta(\tau(\omega))
+\delta(\tau(b))+\delta(\tau(\omega))\\&=&\eta_{A,b}+\delta(\tau(b))+\omega= \eta_b+\omega,
 \end{eqnarray*}
because
\begin{eqnarray*} a-\delta(\tau(a)) - b+\delta(\tau(b))=
\omega  -\delta(\tau(\omega)).  \end{eqnarray*}
Now,  $ \eta_a=\eta_b$ implies $\omega=0$ and hence $a=b$.

For the surjectivity of $\Psi$, as $S$ is  affine, we may assume
that the homomorphism $\rho\colon {\mathbf{G}}^\natural\to G$ is surjective on
$S$\nobd-valued points. Consider then a pair
$({\rho(a)},\nabla_\eta)\in {\cal E}(S)$ with $a$ an $S$\nobd-valued
point of ${\mathbf{G}}^\natural$ and $\eta$ a normal invariant differential of ${\cal
P}_{\rho(a)}$. We defined $\Psi(a)= ({\rho(a)},\eta_a)\in {\cal
E}(S)$. Now, as both $\eta$ and $\eta_a$ are normal invariant
differentials of ${\cal P}_{\rho(a)}$ (i.e.\  they restrict to $dz/z$
on $\mgr {S}$) the differential $\eta-\eta_a$ equals $\omega$ for a
suitable $\omega\in {\boldsymbol{\omega}}_{G'}(S)$. Define $b:=a+\omega$. It
holds $\rho(a)=\rho(b)$ and $\eta_b=\eta_a+\omega=\eta$.

It is also immediate to check that the homomorphism $\Psi$ does not
depend on the choice of the section $\delta$. \qed 

 The definition of Deligne's pairing for the de Rham realizations of
$1$\nobd-motives over a field uses the fact that the pull-back of a biextension
of $1$-motives $(M_1,M_2)$ by $\mgr S$ to the universal extensions
$(M^\natural_1,M^\natural_2)$ admits a  canonical  
 $\natural$\nobd-structure. 
The case of Poincar\'e biextensions can be deduced from
Theorem~\ref{thm.nat}.
However, we prove it separately because we will use in the next section the
explicit description of the canonical $\natural$-structure contained in the proof.

\begin{proposition}\label{pro.natpoincare}
Let $\cal P^\natural$ be the pull-back to $(M^\natural,M^{\prime \natural})$ of the
Poincar\'e biextension $\cal P$  of $(M, M')$ by $\mgr S$. It admits
a   canonical   $\natural$\nobd-structure,  i.e.\  there is a canonical 
connection on the underlying torsor that makes $\cal P^\natural$ a
$\natural$\nobd-biextension of $(M^\natural, M^{\prime \natural})$  by $\mgr S$.
 This is the unique $\natural$-structure on $\cal P^\natural$ if 
$\homo{G^{\natural}}{\mathbb G_a}=0=\homo{G^{\prime\natural}}{\mathbb G_a}$.
 
\end{proposition}
\proof  Denote by ${\cal P}_{\rho} $ the pull-back of ${\cal
P}$ to $(M^\natural, M')$ as well its associated $\mgr S$\nobd-torsor on
${\mathbf{G}}^\natural\times G'$. 
By Proposition~\ref{pro.repE} the identity map on
${\mathbf{G}}^\natural$ provides a $\natural$\nobd-structure $\nabla_2$ on ${\cal
P}_{\rho}$ (viewed as extension of $M'$ by the multiplicative group
over ${\mathbf{G}}^\natural$). 
To check the horizontality condition on $\nu_1$ (see
Definition~\ref{def.12structure}) one uses the isomorphism $\Psi$ in
the proof of  Proposition~\ref{pro.repE}. Indeed, the pull-back via
\begin{eqnarray*}p_{13}\colon {\mathbf{G}}^\natural\times
{\mathbf{G}}^\natural\times G'\to {\mathbf{G}}^\natural\times G' 
\quad (\text{resp.\  }p_{23}, \text{
resp.\  } \mu_{{\mathbf{G}}^\natural}\times \id_{G'})
\end{eqnarray*}
of $({\cal P}_{\rho},\nabla_2)$ is the image via $\Psi$ of the
${\mathbf{G}}^\natural\times {\mathbf{G}}^\natural$-valued point $\rho\circ
p_1$ of $G$ (resp.\ $\rho\circ p_2$, resp.\  
$\rho\circ \mu_{{\mathbf{G}}^\natural}$) and it holds
$p_1+p_2=\mu_{{\mathbf{G}}^\natural}$. Changing the role of $M$ and $M'$ we get a
$\natural$\nobd-$1$\nobd-structure of ${\cal P}^\natural$.
 
To show the uniqueness result it is sufficient to show that any 
$\natural$\nobd-structure $\nabla$ on the trivial biextension of 
$(M^\natural, M^{\prime \natural})$ by $\mgr S$ is trivial 
(cf. \cite{De}, 10.2.7.4.). We are considering the trivial $\mgr
S$\nobd-torsor on ${\mathbf{G}}^\natural\times {\mathbf{G}}^{\prime \natural}$ 
with a connection $\nabla$ such that the morphisms $\nu_1,\nu_2$ in Definition
(\ref{def.natbiext}) are horizontal and the pull-back of $\nabla$ to
${\mathbf{G}}^\natural\times X'$ is trivial as well as the pull-back to $X\times
{\mathbf{G}}^{\prime \natural}$. The connection $\nabla$ is determined by giving a 
global differential $\omega=\omega_1+\omega_2$ on 
${\mathbf{G}}^\natural\times {{\mathbf{G}}^{\prime \natural}}$ where the 
$\omega_i$ depends on the $\natural$\nobd- $i$\nobd-structure associated to 
$\nabla$.
Recall now that $\omega_1$ has to be a global invariant differential  on 
${\mathbf{G}}^\natural_{{\mathbf{G}}^{\prime \natural}}$. We may work Zariski 
locally on $S$ and then assume that the sheaf of differential forms of 
${\mathbf{G}}^\natural$ over $S$ is free. We can then write 
$\omega_1=\sum_j F_j\omega_{1j}$ with $\{\omega_{1,j}\}_j$ (the pull-back) of a
free basis of invariant differentials of  ${\mathbf{G}}^\natural$  and $F_j$  (the
pull-back) of a global section of ${\mathbf{G}}^{\prime \natural}$. The condition on 
$\nu_2$ requires that $F_j$ is additive, i.e.\  it corresponds to a
homomorphism ${\mathbf{G}}^{\prime \natural}\to \agr S$. However 
${\mathbf{G}}^{\prime \natural}$ is extension of $X'\otimes \agr S$ by 
$G^{\prime \natural}$ and then by hypothesis  $F_j$ comes from an additive global 
section of the vector group $X'\otimes \agr S$. 
It is clear that the pull-back of 
$\omega_1$ to $X\times {\mathbf{G}}^{\prime \natural}$ is trivial because $X$ is 
\'etale. Moreover, the condition that the pull-back of $\omega_1$ to 
${\mathbf{G}}^\natural\times X'$ has to be trivial implies that $F_j=0$. 
Hence $\omega_1=0$. In the same way one sees that $\omega_2=0$.
\qed

More generally: 

\begin{theorem}\label{thm.nat}
Let $M_i=[u_i\colon X_i\to G_i], i=1,2$, be two $1$\nobd-motives,
$\cal P$ a biextension of $(M_1,M_2)$ by $\mgr S$ and $\cal
P^\natural$ its pull-back to $(M^\natural_1, M^\natural_2)$. Then $\cal
P^\natural$ admits a   canonical   $\natural$\nobd-structure.   This is the unique 
$\natural$-structure on $\cal P^\natural$ if $\homo{G_i^\natural}{\mathbb G_a}=0$.
 \end{theorem}
\proof This is essentially Deligne's proof in \cite{De}, 10.2.7.4. The
  uniqueness result can be proved as in the previous Proposition. For the existence, observe that \cite{SGA7} VIII, 3.5, implies
that the pull-back homomorphism
\begin{eqnarray*}
\biext{G_1}{[X_2\to A_2]}{\mgr S}\longrightarrow \biext{G_1}{M_2}{\mgr S}
\end{eqnarray*}
is indeed an isomorphism. Hence  $\cal P$ is the pull-back of a
biextension $\tilde{\cal P}$ of $(G_1,[X_2\to A_2])$ by $\mgr S$.
Moreover, $\tilde{\cal P}$ provides a homomorphism
\begin{eqnarray*}\psi\colon G_1\to \uext{[X_2\to A_2]}{\mgr S}=G_2'
\end{eqnarray*}
(cf. \cite{SGA7}, VIII 1.1.4)
and  $\tilde{\cal P}$ is the pull-back via $ \psi\times \id$ of the
Poincar\'e biextension of $(G_2',[X_2\to A_2])$. We define now an
$S$\nobd-group scheme
\begin{eqnarray*}
C:= {\mathbf{G}}^{\prime \natural}_2\times_{G_2'} G_1
\end{eqnarray*}
via the usual homomorphism ${\mathbf{G}}^{\prime\natural}_2\to {G_2'}$ and $\psi$.
The group $C$ is extension of $G_1$ by ${\boldsymbol{\omega}}_{G_2}$. 
Using  Proposition~\ref{pro.repE} one shows that 
\begin{equation}
\label{eq.C} C(S') = \left\{
 \begin{array}{ll}
 ( {g},\nabla ),& g\in G_1(S'), \nabla \text{ a }
\natural\text{-structure on the corresponding} \\
&\text{extension  } {\cal P}_{g } \text{ of } M_2 \text{ by }\mgr
{S'}
\end{array}\right\}.
\end{equation}
 Define now a homomorphism $u_C\colon X_1\to C,$  as $u_C(x)=(u_1(x), \nabla^0)$
where  $\nabla^0$ denotes the trivial connection on  ${\cal
P}_{u_1(x)}$. Observe that, by definition of biextensions of
complexes, the pull-back of ${\cal P}$ to $X_1\times G_2$ is
isomorphic to the trivial biextension. In this way $u_C\colon X_1\to
C$ becomes an extension of $M_1$ by the vector group ${\bm
\omega}_{G_2}$. Using the universal property of the universal
extension $M_1^\natural=[X_1\to {\mathbf{G}}^\natural_1]$ of 
$M_1$, $u_C \colon X_1\to C$ is the push-out of $M_1^\natural$ for a suitable
homomorphism ${\boldsymbol{\omega}}_{G_1'}\to
{\boldsymbol{\omega}}_{G_2}$. Denote by $\Gamma$ the induced homomorphism 
${\mathbf{G}}^\natural_1\to C$. It is clear
that the image via $\Gamma$ of the identity of  ${\mathbf{G}}^\natural_1$
provides a ${\mathbf{G}}^\natural_1$\nobd-valued point of $C$ that
corresponds, because of (\ref{eq.C}), to a
$\natural$\nobd-$2$\nobd-structure on  the pull-back of ${\cal P}$
to $({\mathbf{G}}^\natural_1,M_2)$ and hence on ${\cal P}^\natural$. In a
similar way, one gets a $\natural$\nobd-$1$\nobd-structure on ${\cal
P}^\natural$ and hence the canonical $\natural$\nobd-structure we are looking
for. \qed

\begin{remark}\label{rem.notunique}
The uniqueness result in \cite{De}, 10.2.7.4 
(see also the proof of Propositions~\ref{pro.natpoincare})
depends on the fact that $\homo{G^\natural}{\mathbb G_a}=0$. This is not true in 
general. Indeed $\homo{G^\natural}{\mathbb G_a}$
is the kernel of the push-out homomorphism $\homo{\mathbf \omega_{A'} }{\mathbb G_a}
\to \ext{G}{\mathbb G_a}$. This map is an epimorphism because of the universal property
of universal extensions (cf. (\ref{def.epsilon})). It can not be an isomorphism
when $\homo{\mathbf \omega_{A'} }{\mathbb G_a}$ is bigger than
${\mathrm {Hom}}_{\cal O_S}(\mathbf \omega_{A'} ,\mathbb G_a)$. As an example,
over a field $k$ of characteristic $p>0$, for $\mathbf \omega_{A'}=\mathbb G_a$, 
the homomorphisms of $k$-group schemes $\mathbb G_a\to \mathbb G_a$ correspond 
to polinomials of the type $\sum_i a_i x^{p^i}, a_i\in k$, while the
homomorphisms of vector groups $\mathbb G_a\to \mathbb G_a$ correspond to 
linear polinomials $ax$, $a\in k$.
\end{remark}

\section{Deligne's pairing $\Phi$}\label{sec.perf}
Let $M_1, M_2$ be two $S$\nobd-$1$\nobd-motives, ${\cal P}$  a biextension of
$(M_1, M_2)$ by $\mgr S$ and ${\cal P}^\natural$ the pull-back of
${\cal P}$ as biextension of  $(M_1^\natural,M_2^\natural)$ by $\mgr S$.
Following Deligne, denote $\uLie{{\mathbf{G}}^\natural_i}$ by $\TT_\dR(M_i)$.
We know from Theorem \ref{thm.nat} that ${\cal P}^\natural$ admits a
canonical $\natural$\nobd-structure. Hence ${\cal P}^\natural$ is
equipped with a canonical connection $\nabla$. Consider now the
curvature form of $\nabla$ (see, for example, \cite{MM} I, 3.1.4).
It is an invariant $2$\nobd-form on ${\mathbf{G}}^\natural_1 \times
{\mathbf{G}}^\natural_2 $; hence it gives an alternating pairing $R$ on
\begin{eqnarray*}
\uLie{{\mathbf{G}}^\natural_1\times {\mathbf{G}}^\natural_2}=
\uLie{{\mathbf{G}}^\natural_1}\oplus \uLie{{\mathbf{G}}^\natural_2}=
\TT_\dR(M_1)\oplus \TT_\dR(M_2)\end{eqnarray*}
 with values in $\uLie{\mgr S}$. As the restrictions of $R$ to
$\uLie{{\mathbf{G}}^\natural_i}$, $i=1,2$,  are trivial it holds
\begin{eqnarray*}
R(g_1+g_2,g_1'+g_2')=\Phi(g_1,g_2')-\Phi(g_2,g_1')
\end{eqnarray*}
with $\Phi\colon \TT_\dR(M_1)\otimes  \TT_\dR(M_2) \to \uLie{\mgr S}$
a bilinear map. We will show in this section that  $\Phi$ is a
non-degenerate pairing when $M_1, M_2$ are Cartier duals and ${\cal
P}$ is the Poincar\'e biextension.  This result has been proved by
Deligne for $S=\spec \CC$  and by Coleman for  abelian
schemes and $S$ flat over $\ZZ$. See also \cite{Gr}, V \S4. Both proofs 
are based on the comparison with another perfect pairing and do not work
in our general case. 

Let in the following  $M=[u\colon X\to G], M'=[u'\colon X'\to G']$
be Cartier duals, $\cal P$ the Poincar\'e biextension of $(M,M')$ and
\begin{eqnarray}\label{pairing2}
\Phi\colon \TT_\dR(M)\otimes  \TT_\dR(M') \to \uLie{\mgr S}
\end{eqnarray}
Deligne's pairing. Recall that we have vectorial extensions of $M$
and $M'$
\begin{eqnarray}\label{seq.wEMG}
\xymatrix @C=17pt{
 &0\ar[r]&{\boldsymbol{\omega}}_{G'}\ar[r]^(.5)i & M^\natural \ar[r]^(.5)\rho& M\ar[r]&
 0, & 
 0\ar[r]&{\boldsymbol{\omega}}_{G}\ar[r]^(.5){i'} & M^{\prime \natural} \ar[r]^(.5){\rho'}& M'\ar[r]&
 0,
}
\end{eqnarray}
 with $M^\natural=[X\to {\mathbf{G}}^\natural ]$, $M^{\prime \natural}=[X'\to {\mathbf{G}}^{\prime \natural}]$ the
 universal extensions of $M$, $M'$. Recall that ${\cal P}^{\natural}$
 denotes the  pull-back of ${\cal P}$ to
 $(M^\natural,M^{\prime \natural})$ endowed with its canonical
 $\natural$-structure. 
We showed in Proposition~\ref{pro.natpoincare} that  ${\cal P}^{\natural}$ is
the sum of  $({\cal P}_{\rho},\nabla_2)$ 
and $({\cal P}_{\rho'},\nabla_1)$ (after suitable pull-backs) where  $({\cal
  P}_{\rho},\nabla_2)$ is  the $\natural$\nobd-extension of $M'$ by the
multiplicative group over ${\mathbf{G}}^\natural$ that corresponds to the
identity map on ${\mathbf{G}}^\natural$ via the isomorphism $\Psi$ in
(\ref{def.psi}). Similarly for  $({\cal P}_{\rho'},\nabla_1)$.

\begin{lemma}\label{lem.alpha} Let
$\alpha_{G'}$ be the invariant differential of $G'$ over ${\bm
\omega}_{G'}$ that corresponds to the identity map on ${\bm
\omega}_{G'}$. The restriction of $({\cal P}_{\rho },\nabla_2)$ to
${\boldsymbol{\omega}}_{G'}$ via $i\colon {\boldsymbol{\omega}}_{G'}\to {\mathbf{G}}^\natural$ in
(\ref{seq.wEMG}) is isomorphic to the trivial extension of $ M'$ by
the multiplicative group over ${\boldsymbol{\omega}}_{G'}$ equipped with the
connection associated to $\alpha_{G'}$.
\end{lemma}
\proof (See also \cite{Co2}, Lemma 2.0 for the case $M=[0\to A]$.)
Recall that we have the following arrows
\begin{eqnarray*}
\xymatrix{ {\mathbf{G}}^\natural({\mathbf{G}}^\natural)\ar[r]^(.5)F& {\mathbf{G}}^\natural({\bm
\omega}_{G'})&\ar[l]_{ H } {\boldsymbol{\omega}}_{G'}({\boldsymbol{\omega}}_{G'}) }
\end{eqnarray*}
where the $F(f)=f\circ i$ and $H(h)=i\circ h$. In terms of
$\natural$\nobd-extensions of $M'$ by the multiplicative group, the
homomorphism $F$ is the base-change via $i$, while $H$ associates to
a differential $\eta$ the trivial extension of $M'$ by the
multiplicative group over ${\boldsymbol{\omega}}_{G'}$ endowed with the
connection associated to $\eta$. As $f(\id)=i=H(\id)$, the
restriction of $({\cal P}_{\rho},\nabla_2)$ to ${\boldsymbol{\omega}}_{G'}$
is isomorphic to the trivial extension of $M'$ by the multiplicative
group over ${\boldsymbol{\omega}}_{G'}$ equipped with the connection
associated to $\alpha_{G'}$. 
\qed

Changing the role of $M$ and $M'$, denote by $\alpha_{G}$
 the invariant differential of $G$ over ${\boldsymbol{\omega}}_{G}$
that corresponds to the identity map on ${\boldsymbol{\omega}}_{G}$. The
restriction  $({\cal P}_{\rho'},\nabla_1)$ to ${\boldsymbol{\omega}}_{G}$ is
isomorphic to the trivial extension of $M'$ by the multiplicative
group over ${\boldsymbol{\omega}}_{G}$ equipped with the connection
associated to $\alpha_{G}$.

In order to study the curvature forms of the connections $\nabla_i$ we start
considering the curvatures of $\alpha_G$ and $\alpha_{G'}$.
We will use in the following the same notation for a locally
free sheaf and its associated vector group.

\begin{lemma}\label{lem.dalpha}
The curvature of $\alpha_G$ provides a perfect pairing
\begin{eqnarray*}
d\alpha_G\colon~{\boldsymbol{\omega}}_{G}\otimes \uLie{G} \longrightarrow
\uLie{\mgr S}
\end{eqnarray*}
that is the usual duality.
\end{lemma}
\proof We may work locally.  Let
$\omega_1,\dots,\omega_g$ be a basis of invariant differentials of
$G$. We have ${\boldsymbol{\omega}}_{G}=\spec{{\cal O}_S[x_1,\dots,x_n]}$
where $x_i$ is the basis of $\uLie{G}$ dual to $(\omega_i)_i$. An
$S'$-valued point of ${\boldsymbol{\omega}}_{G}$ corresponds to a $g$-tuple
$(a_i)_i\in \Gamma({\cal O}_{S'},S')^g$, hence to the invariant
differential $\sum_ia_i\omega_i$ of $G$ over $S'$. Therefore
$\alpha_G=\sum_ix_i\omega_i$. In particular, its curvature form
$d\alpha_G=\sum_idx_i\wedge \omega_i$ provides a pairing
\[\uLie{{\boldsymbol{\omega}}_{G} }\otimes \uLie{G} \to {\cal O}_S
\]
that is the usual duality, once identified $\uLie{{\boldsymbol{\omega}}_{G} }$ with
${\boldsymbol{\omega}}_{G}$.\qed

\begin{theorem}\label{thm.perfM}
Let $M$ be an $S$-$1$-motive. Then Deligne's pairing in $(\ref{pairing2})$ is perfect.
\end{theorem}
\proof Observe that the biextension ${\cal P}_{\rho'}$ can be defined also as the
pull-back of the Poincar\'e biextension of $([X\to A],G')$ to
$(M,{\mathbf{G}}^{\prime \natural})$ and the biextension ${\cal P}^{\natural}$ is the
pull-back of ${\cal P}_{\rho'}$ via $$(\rho\times\id)\colon M^\natural
\times {\mathbf{G}}^{\prime \natural}\to M \times {\mathbf{G}}^{\prime \natural}$$ together with a
suitable trivialization on $ {\mathbf{G}}^\natural\times X'$. Furthermore, we
have an exact sequence of ${\mathbf{G}}^{\prime \natural}$-group schemes
\begin{eqnarray*}
\xymatrix{
 0\ar[r]& {\boldsymbol{\omega}}_{G'} \times_ S{{\mathbf{G}}^{\prime \natural} }\ar[r] & {\cal
 P}^{\natural}_{}
\ar[r] & {\cal P}_{\rho'}\ar[r]& 0. }
\end{eqnarray*}
This assures that, after pull-back to $ {\boldsymbol{\omega}}_{G'} \times_
S{{\mathbf{G}}^{\prime \natural} } $, the $\natural$-$1$-structure of ${\cal
P}^{\natural}$  is the trivial connection because it
 comes from the connection $\nabla_1$ on  ${\cal P}_{\rho'}$.
 Lemma~\ref{lem.alpha} implies that  after pull-back to
${\boldsymbol{\omega}}_{G'}\times {\mathbf{G}}^{\prime \natural}$, the $\natural$-$2$-structure
of ${\cal P}^{\natural}$ is the connection associated to
 the invariant differential $\rho^{\prime *}\alpha_{G'}$ of ${\mathbf{G}}^{\prime \natural}$.
Hence, the restriction of the curvature form of $\nabla$ to
\[\uLie{ {\boldsymbol{\omega}}_{G'} }\otimes \uLie{ {\mathbf{G}}^{\prime \natural} }=
{\boldsymbol{\omega}}_{G'}\otimes \uLie{ {\mathbf{G}}^{\prime \natural} } \] is $d(\rho^{\prime
*}\alpha_{G'})$. This says that the homomorphisms $\iota$ and $g'$
in the following sequences of Lie algebras deduced from
(\ref{seq.wEMG})
\begin{eqnarray*}
\xymatrix{
 0\ar[r]&{\boldsymbol{\omega}}_{G'}\ar[r]^(.5){\iota}&\uLie{{\mathbf{G}}^\natural}\ar[r]^(.5){g}&\uLie{G}\ar[r]&
 0,\\
 0& \ar[l] \uLie{G'}&\uLie{{\mathbf{G}}^{\prime \natural}}\ar[l]^{g'} & {\boldsymbol{\omega}}_{G}\ar[l]^{\iota'}& 0\ar[l],
}
\end{eqnarray*}
are transposes of each other with respect to the pairing
$d\alpha_{G'}$  and  $\Phi$. Changing the role of $M$ and $M'$,
we get that $\iota'$ and $g$ are transposes of each other with
respect to  $d\alpha_G$ and  $\Phi$. The perfectness of
$d\alpha_{G}$ and $d\alpha_{G'}$ was proved in
Lemma~\ref{lem.dalpha}. Hence also $\Phi$ is perfect. \qed

\begin{example}\label{ex.deligne} Case $A=A'=0$ and $T,T'$ split.
In this case it is possible to give an explicit description of Deligne's pairing.
 Let $M$ be of the form $
[u\colon \ZZ^r\to \mgr S^d]$. Then $M'$ is of the form $[u'\colon
\ZZ^d\to \mgr S^r]$ and
\begin{eqnarray*}
M^\natural =[(\iota,u)\colon \ZZ^r\to \agr S^r\times \mgr S^d], \quad
M^{\prime \natural} =[(\iota',u')\colon \ZZ^d\to \agr S^d\times\mgr S^r],
\end{eqnarray*} where we write $\agr S^i$ in place of ${\bm
\omega}_{\mgr S^i}$, for $i=r,d$, and $\iota$ (resp.\  $\iota'$) sends an
$r$-tuple (resp.\  a $d$-tuple) $\bm n$ to ${\bm n}$.

Suppose  $r=1, d=0$. The pull-back of $\cal P$ to the universal
extensions $(M^\natural,M^{\prime \natural})$ is the trivial biextension of $(\agr
S, \mgr S)$ by ${ \mgr S}$ together with two trivializations
\begin{eqnarray*}
\tau_1\colon \quad \ZZ\times \mgr S\to {\cal P}^\natural= \agr
S\times \mgr S\times \mgr S,
\quad (n,b)\mapsto (n, b, b^n),\nonumber \\
\tau_2 \colon \quad \agr S\times 0\to {\cal P}^\natural= \agr
S\times \mgr S\times \mgr S, \quad (a,0)\mapsto (a, 1, 1).
\end{eqnarray*}
that coincide on $\ZZ\times 0$. Observe that the above biextension
of complexes is \emph{not} trivial. In order to describe the canonical
$\natural$\nobd-structure on ${\cal P}^\natural$, we start constructing the
global differential $\omega$ on $\agr S\times \mgr S$ associated to
a connection on ${\cal P}^\natural$. Let $\agr S=\spec
{{\cal O}_S[x]}, \mgr S=\spec {{\cal O}_S[t, 1/t]}$ and $dx, dt/t$
be  the usual  invariant differentials; it will be
\begin{eqnarray*}\omega=f(x,t)dx+g(x,t)dt/t \quad \text{with}
 \quad f,g\in \Gamma(S, {\cal O}_S)[x,t, 1/t].
\end{eqnarray*}
Recall that  $\nu_1, \nu_2$ in Definition~\ref{def.natbiext} are
asked to be horizontal. An easy computation shows that necessarily
$f=0$ and   $g$ is ``additive'' in $x$ and does not depend on $t$. Observe that 
over a filed of characteristic zero $g(x)=ax$, while 
over a field of characteristic $p$, we could a priori have the case $g(x)=x^p$.  
Hence
\begin{eqnarray*}\omega=g(x) dt/t , \quad   g(x)\in  \Gamma(S, {\cal O}_S)[x].
\end{eqnarray*}
We use now the hypothesis that the pull-back of ${\cal P}^\natural$
to $\ZZ\times \mgr S$ is isomorphic to the trivial
$\natural$\nobd-biextension. The trivialization $\tau_1$ restricted
to $\{n\}\times \mgr S$ is $(n,b)\to (n,b,b^n)$ and ${\cal
P}_{|\{n\}\times \mgr S}^\natural$ has to be isomorphic to
the trivial $\natural$\nobd-extension of $\mgr S$ by itself. The
pull-back of $\omega$ to $\mgr S$ is $g(n)dt/t$  and this has to equal 
$ndt/t$   (see Example~\ref{ex.1}).  
Therefore  $g(x)=x$, the canonical
$\natural$\nobd-structure on ${\cal P}^\natural$ is given by the
differential $xdt/t$ and its curvature by $dx\wedge dt/t$. In
particular, Deligne's pairing
\begin{eqnarray*}
\Phi\colon \uLie{\agr S}\otimes \uLie{\mgr S}\to \uLie{\mgr S}
 \end{eqnarray*}
 is non degenerate.

In the general situation one proceed in a similar way and the global
differential associated to the canonical $\natural$\nobd-structure on
${\cal P}^\natural$ is of the form
\begin{eqnarray*}
\omega=\sum_{i=1}^r x_idz_i/z_i +\sum_{j=1}^d y_jdt_j/t_j
\end{eqnarray*}
where the parameters $x_i$ (resp.\  $y_i$) refer to $\agr S^r $ (resp.\ 
to $\agr S^d$) and the $z_i$ (resp.\  $t_j$) refer to $\mgr S^r $
(resp.\  to $\mgr S^d$).  Again the curvature gives rise to a non
degenerate pairing whose matrix is built of two blocks, each one
involving a torus and its group of characters.
\end{example}

\begin{remark} As one expects, Deligne's pairing is compatible 
with weight filtration. To see this fact we have to work locally because, 
in general, we have no canonical morphisms 
${\mathbf G}^\natural\to A^\natural$. We assume then that there exist sections 
$\delta$ of $\bar\tau\colon{\boldsymbol{\omega}}_{G'}\to
{\boldsymbol{\omega}}_{T'}$
and  $\delta'$ of $\bar\tau'\colon{\boldsymbol{\omega}}_{G}\to
{\boldsymbol{\omega}}_{T}$ as in the proof of Proposition~\ref{pro.repE}. We get then a morphism 
${\mathbf{G}}^\natural\stackrel{\id-\delta\tau}{\to} G^\natural \to
A^\natural$ and similar for ${\mathbf{G}}^{\prime \natural} $ so that we have  an extension
\begin{eqnarray}\label{seq.TGA}
\xymatrix{
 0\ar[r]&T\times T'\times {\boldsymbol{\omega}}_{T'} \times {\boldsymbol{\omega}}_{T}\ar[r]&
{\mathbf{G}}^\natural\times {\mathbf{G}}^{\prime \natural} \ar[r]& 
A^\natural \times A^{\prime \natural} \ar[r]& 0.
}\end{eqnarray}
Recall now that Deligne's pairing $\Phi$ for $M, M'$ is defined via
the canonical $\natural$-structure $\nabla$ on ${\cal P}^\natural$, 
the pull-back of the Poincar\'e biextension of $(M, M')$ to $(M^\natural,
M^{\prime\natural})$.
Let now ${\cal P}^{A\natural}$ be the pull-back of the Poincar\'e
biextension of $A, A'$ to $A^\natural, A^{'\natural}$ and $\nabla^A$ its
canonical $\natural$-structure.  We denote by $({\cal
  P}^{A\natural},\nabla^A)$ also its pull-back to ${\mathbf G}^\natural\times
{\mathbf G'}^\natural$. With the notations of  Lemma~\ref{lem.alpha}, let 
$\delta^*\alpha_{T'}$ be the invariant
differential of $G'$ over ${\mathbf G}^\natural$ associated to 
$\delta\tau\in {\boldsymbol{\omega}}_{G'}({\mathbf G}^\natural)$ and also its
pull-back to ${\mathbf G}^{\prime \natural}$.
Similarly for $\delta^{\prime   *}\alpha_T$.  Following the constructions in
the proofs of Propositions~\ref{pro.repE},~\ref{pro.natpoincare}, for
$\id_{\mathbf{G}^\natural}=\id_{\mathbf{G}^\natural}-\delta\tau\oplus \delta\tau$, and the ``dual'' one, we get that  
\[ ({\cal P}^\natural, \nabla)= ({\cal P}^{A\natural},\nabla^A)+
(0,\delta^*\alpha_{T'}+\delta^{'*}\alpha_T) \]
as $\natural$-biextension of 
$({\mathbf G}^{\natural}, {\mathbf G}^{\prime \natural})$ by $\mgr S$.
From this decomposition and the sequence of Lie algebras of (\ref{seq.TGA})
we get that Deligne's pairing for $M, M'$ lies in the middle of a diagram 
where on the right we have Deligne's pairing for $A, A'$ and  on the left, 
after re-ordering the summands,  we have Deligne's pairing for $[X\to T],
[X'\to T']$. The latter pairing has locally the concrete description given in 
the previous example. 
\end{remark}

\section{$\natural$-extensions of $1$-motives}\label{sec.natext}

This section contains results on $\natural$\nobd-extensions of a
$1$\nobd-motive by the multiplicative group that generalize what
happens in the classical case of abelian schemes. In particular, we show that
$\natural$\nobd-extensions are no longer sufficient to describe universal
extensions of $1$-motives. 

 Let notations be as in section~\ref{sec.ue}. Recall that $A'$
represents the functor $S'\leadsto  \ext{A_{S'}}{\mgr {S'}}$.

\begin{theorem}[\cite{MM}]\label{thm.extnatA}
The universal extension $A^{\prime \natural}$ of the abelian scheme $A'$ over
$S$ represents the functor ${\cal F}_A\colon S'\leadsto  \extnat {A_{S'}}{\mgr {S'}}$.
\end{theorem}
As a consequence we can interpret $S'$\nobd-valued points of
$A^{\prime \natural}$ as $\natural$\nobd-extensions of $A$ by ${\mathbb G}_m$
over $S'$.

\begin{lemma} \label{lem.extnatMA}
The pre-sheaf for the flat topology
$
S'\leadsto  \extnat {[X_{S'}\to A_{S'}]}{\mgr {S'}}
$
is a sheaf represented by ${G}^{\prime \natural}$, the universal extension of
$G'$.
\end{lemma}
\proof It is known (cf. \cite{ABV}, 2.2.1) that
$ G^{\prime \natural}=A^{\prime \natural}\times_{A'} G'$. In particular
\begin{eqnarray*}{G}^{\prime \natural}(S')=\{(x,y)\in A^{\prime \natural}(S')\times G'(S') ~~
\text{inducing the same $S'$\nobd-valued point on }
A'\}.\end{eqnarray*} One checks immediately that $\extnat
{[X_{S'}\to A_{S'}]}{\mgr {S'}}$ is the group
\begin{equation*}
\left\{
 \begin{array}{ll}(x,y)\in &
\extnat{A_{S'}}{\mgr {S'}}\times \ext{[X_{S'}\to A_{S'}]}{\mgr {S'}}\\
&\text{inducing the same extension of }A_{S'}\text{ by }\mgr {S'}
\end{array}\right\}.
\end{equation*}
Recalling that $G'$ represents the functor 
$S'\leadsto  \ext{[X_{S'}\!\to \!A_{S'}]}{\mgr {S'}},
$
the conclusion follows from Thm.~\ref{thm.extnatA}. \qed

For $G$ a semi-abelian scheme over $S$ with maximal subtorus $T$, it
is no longer true that the pre-sheaf for the flat topology
$S'\leadsto  \ext{G_{S'}}{\mgr {S'}}$ is a sheaf. Indeed, for
$G=T$ the associated sheaf is trivial. However the functor
\begin{eqnarray*}{\cal F}_G \colon S'\leadsto  \extnat{G_{S'}}{\mgr
{S'}}
\end{eqnarray*}
is still a sheaf   if we restrict to a suitable site.

\begin{lemma}\label{lem.iso}
Suppose $S$ flat over $ \ZZ$. Let $(E_1,\nabla_1), (E_2,\nabla_2)$
be $\natural$\nobd-extensions of the semi-abelian scheme $G$ by
$\mgr S$. Suppose given two horizontal isomorphisms of extensions
$f,g\colon E_1\to E_2$. Then $f=g$.
\end{lemma}
 \proof The result is trivially true if $G=A$
because $g^{-1}f$ (resp.\  $f^{-1}g$) is an automorphism of the
extension $E_1$ (resp.\  $E_2$) and hence it coincides with the
identity map.

Suppose now that the abelian part $A$ is trivial and $G=T$. We may
work (fppf) locally on $S$ an then suppose that both $E_i$ are the
trivial extension $E^0=\mgr S\times_S T$. The isomorphisms $f,g$
correspond, respectively, to characters $a,b\colon T\to \mgr S$. Let
$dz/z+\omega_i$ be the normal invariant differential on $E^0$
associated to the connection $\nabla_i$, $i=1,2$, where $\omega_i$
are invariant differentials of $T$. The horizontality condition says
that $f^*\nabla_2=g^*\nabla_2=\nabla_1$. Hence
\begin{eqnarray*}
dz/z+da/a+\omega_2=dz/z+db/b+\omega_2=dz/z+\omega_1, \quad{\rm where~} d
a/a:=a^*(dz/z).
\end{eqnarray*} 
Therefore $da/a=db/b$ and this implies $a=b$
because of the hypothesis on $S$.\footnote{This proof does not work
in positive characteristic $p$.  For example, given characters $a,b$
of $T$ one has $da/a=d(ab^p)/ab^p$.}

In the general situation let $E_{iT}$ be the pull-back of $E_i$ via
$T\to G$. It is clear that we have exact sequences
\begin{eqnarray*}
\xymatrix{0\ar[r] & E_{iT}\ar[r] &E_i\ar[r] & A\ar[r]&0.}
\end{eqnarray*}
The isomorphisms $f,g$ induce isomorphisms of tori $f_T,g_T\colon
E_{1T}\to E_{2T}$ and $f_T=g_T$ because of what we explained above.
Hence $g^{-1}_Tf_T=\id_{E_{1T}}$ and $g^{-1}f$ is an automorphism of
the extension $E_1$ that necessarily coincides with the identity of
$E_1$ because there exist no non-trivial homomorphisms of $A$ to
$E_{1T}$.
\qed

Let  $S$ be a scheme flat over $\ZZ$, ${\bf Sch}/S$ the category of
$S$-schemes and   ${\bf Fl}/S$ the full 
subcategory of ${\bf Sch}/S$ consisting of those $S$-schemes flat over $\ZZ$.
 Observe that if $S=\spec k$ with $k$ a field of characteristic $0$, all
 $S$-schemes are flat over $\ZZ$, hence, ${\bf Fl}/S$ and ${\bf Sch}/S$
 coincide. More generally,
this is true if the following hypothesis holds:
\medskip

$(*)$\quad \quad {\it All residue fields of $S$ have characteristic $0$.}
 \medskip

 We prove now that ${\cal F}_G$ is a sheaf on the site $({\bf Fl}/S)_{\rm fl}$.

\begin{proposition}\label{pro.extnatG}
Let $G$ be a semi-abelian scheme over $S$ and suppose $S$ flat over
$\ZZ$. Then the functor
\begin{eqnarray*}{\cal F}_G\colon S'\leadsto \extnat{G_{S'}}{\mgr {S'}}
\end{eqnarray*} 
is a sheaf   on $({\bf Fl}/S)_{\rm fl}$. 
\end{proposition}
\proof We start showing that it is a separated pre-sheaf. Let
$(E,\nabla)$ be a $\natural$\nobd-extension of $G$ by the
multiplicative group over $S'$. Suppose it trivializes over a
covering $\{S'_j\}_j$ of $S'$. Hence for any index $j$  we have an
isomorphism $\varphi_j$ of $(E_{S_j'},\nabla)$ with the trivial
$\natural$\nobd-extension $(E^0,\nabla^0)$ of $G_{S_j'}$ by $\mgr
{S_j'}$. Moreover, $\varphi_i^{-1}\varphi_j,
\varphi_j^{-1}\varphi_i$ are horizontal automorphisms of $E$ over
$S_{ij}':=S_i'\times_{S'}S_j'$. By the previous lemma we conclude
that $\varphi_j=\varphi_i$ over $S_{ij}'$ and hence these
isomorphisms descend to a $\varphi\colon E\to \mgr {S'}\times_{S'}G$
and $E$ is isomorphic to the trivial extension of $G_{S'}$ by $\mgr
{S'}$. Now, $\nabla$ corresponds to a global invariant differential
$\omega$ on $G$ over $S'$. By hypothesis $\omega=da_j/a_j$ over
$S_j'$ for a suitable homomorphism $a_j\colon G_{S'_j}\to \mgr
{S'_j}$ with $da_j/a_j=da_i/a_i$ over $S_{ij}'$. Hence the $a_i$
provides a homomorphism $a\colon G_{S'}\to \mgr {S'}$ and
$\omega=da/a$.

 We finish the proof invoking \cite{Mi} II, 1.5. First we
show that ${\cal F}_G$ is a sheaf for the Zariski topology and then
that
\begin{eqnarray*}
\xymatrix{
{\cal F}_G(U)\ar[r] & {\cal F}_G(U')\ar@<2pt>[r]\ar@<-2pt>[r] &
 {\cal F}_G(U'\times_U U')
}
\end{eqnarray*}
is exact with $U, U'$ both affine and $U'$ flat over $U$.

 Let $\{S_j'\}_j$ be a Zariski\nobd-covering of $S'$ and let
$(E_j,\nabla_j)$ be $\natural$\nobd-extensions over $S_j'$ with
isomorphisms $\varphi_{ij}$ between $(E_i,\nabla_i)$ and
$(E_j,\nabla_j)$ over $S_{ij}'$. Thanks to Lemma~\ref{lem.iso} the
$\varphi_{ij}$ satisfy the usual cocycle condition and hence the
$E_j$ glue together providing an extension $E$ of $G_{S'}$ by $\mgr
{S'}$. The $\natural$\nobd-structure can be defined locally on $S'$
and hence we are done.

Suppose now $U$ affine and let $U'$ be an affine scheme faithfully
flat and locally of finite type over $U$. Let  $(E_{U'},\nabla')$ be
a $\natural$\nobd-extension over $U'$ that provides isomorphic
$\natural$\nobd-extensions on $U'\times_{U}U'$ via the projection
morphisms. Again, the cocycle condition is satisfied because of
Lemma~\ref{lem.iso} and the effectiveness of descent data in the
affine case permits to conclude that $E_{U'}$ descends to an
extension $E$ of $G_{U}$ by $\mgr {U}$. Because of the affine
hypothesis, $E$ admits a $\natural$\nobd-structure. Hence we are
reduced to see that the $\natural$\nobd-structure descends in the
case when $E$ is the trivial extension. But this is obvious because
${\bm\omega}_{G}$ is a sheaf. \qed

It is an easy consequence of the above proposition that

\begin{corollary}\label{cor.extnatM}
Let $S$ be a scheme flat over $\ZZ$. Then  ${\cal F}_M\colon S'\leadsto  \extnat {M_{S'}}{\mgr {S'}}$
  is a sheaf on   on $({\bf Fl}/S)_{\rm fl}$.  \end{corollary}

As we have already remarked, the functor
\[({\bf Sch}/S)^0\longrightarrow ({\bf Sets}), \quad  
S'\mapsto \extnat{M_S'}{\mgr {S'}}
\] 
is not, in general, a sheaf for the flat topology. 
Let denote by $\uextnat{M}{\mgr S}$  the associated sheaf. Its restriction
to  $({\bf Fl}/S)_{\rm fl}$ is the sheaf ${\cal F}_M$ in
Corollary~\ref{cor.extnatM}.

Let  ${\cal H}(M):=\uhomo{M}{\mgr S}$ and ${\cal H}^\nabla(M):=\uhomon{M}{\mgr
  S}=\ker({\cal H}(M)\to {\bm \omega}_G )$.
The sheaf $\uextnat{M}{\mgr S}$  fits in the following exact sequence
 \begin{eqnarray}\label{seq.extnatM}
\xymatrix{
0\ar[r]&{\cal H}^\nabla(M) \ar[r] & {\cal H}(M)\ar[r]& {\bm
\omega}_G\ar[r]^(0.3){j} & \uextnat{M}{\mgr S}\ar@{->>}[r] & \uext{M}{\mgr S}  }
\end{eqnarray}
that generalizes the one in (\ref{seq.homextnatext}) and the one in \cite{MM},
II.4.2. 
 The exactness on the left is assured by definition of the first
  sheaf,   while
the map on the right is an epimorphism because of the commutativity of the 
following diagram
\begin{eqnarray}\label{dia.alphagamma}
\xymatrix @C=12pt{
 0\ar[r]&{\boldsymbol{\omega}}_{A}\ar[r]\ar[d]&G^{\prime \natural}=
\uextnat{[X\to A]}{\mgr S}\ar[r]\ar[d]^\alpha &
G'=\uext{[X\to A]}{\mgr S}\ar[r]\ar@{->>}[d]^\gamma &0\\
 & {\boldsymbol{\omega}}_{G} \ar[r]^(.4){\bar\iota}& \uextnat{M}{\mgr S}
 \ar[r] & \uext{M}{\mgr S} &  }
\end{eqnarray}
where the upper sequence is the one describing $G^{\prime \natural}$ as universal 
extension of $G'$ by $\VV(G')=  {\boldsymbol{\omega}}_{A}$  and
$\gamma$, $\alpha$ are the pull-back homomorphisms.

The remaining part of this section is devoted to prove the following
result:

\begin{proposition}\label{pro.extnatex}
Let $M=[u\colon X\to G]$ be an $S$\nobd-$1$\nobd-motive. The
sequence
\begin{eqnarray*}
\xymatrix{0\ar[r]& 
\displaystyle\frac{ {\cal H}^\nabla(T) }{ {\cal H}^\nabla(T) \cap {\cal H}(G) } 
\ar[r]^{\bar\delta}&\uextnat A{\mgr S}\ar[r]^(.5)\alpha & \uextnat G{\mgr
S}\ar[r]^(.5)\beta& \uextnat T{\mgr S} \ar[r]&0}
\end{eqnarray*}
is exact, where $\alpha,\beta$ are the usual pull-back
homomorphisms,  ${\cal H}^\nabla(T) \cap {\cal H}(G)$ denotes the pull-back of
${\cal H}^\nabla(T)$ to ${\cal H}(G)$ via the monomorphism $\uhomo{G}{\mgr
  S}\to\uhomo{T}{\mgr S}$.  
\end{proposition}
  The morphism $\bar\delta$ above is induced by 
\[\delta\colon \uhomon{T}{\mgr S}\to \uextnat A{\mgr S}, \quad 
x\mapsto [(G_x,\nabla_x)]\]
where $G_x$ is the extension obtained as push-out of $G$ with
respect to the character $x\colon T\to \mgr S$.
The connection $\nabla_x$ is induced by the canonical invariant
differential of $G_x$ that pulls back to $dz/z$ on $\mgr S$ and to $0$ on $G$.

To characterize the kernel of $\beta$ we will need the
following result:

\begin{lemma}\label{lem.normaldiff} Denote by $f_x\colon G\to G_x$
the push-out of $G$ with respect to a character $x\colon T\to \mgr
S$. Define the homomorphism $\sigma_G\colon G\to G_x\times_A G$
via   $f_x$ and the identity on $G$.
\begin{eqnarray*}
\xymatrix{ 0\ar[r] & T\ar[r]^(.5)j \ar[d]^x & G\ar[r] \ar[d]^{f_x} & A
\ar@{=}[d]
\ar[r] & 0 \\
0\ar[r] & \mgr S\ar[r]^(.5)\iota \ar@{=}[d] & G_x\ar[r]^(.5)g & A  \ar[r] & 0 \\
 0\ar[r] & \mgr S\ar[r]^(.5)\iota & G_x\times_A G\ar[r]\ar[u]^{p_{G_x}} &
 G  \ar@{.>}@/_/[l]_{\sigma_G}\ar[u]\ar[r] & 0.}
\end{eqnarray*}
\begin{itemize}
\item Given a normal invariant differential $\eta_x$ on $G_x$ it
holds $j^*f_x^*(\eta_x)=dx/x$.
 \item Given  a normal invariant
differential $\eta$ of $G_x\times_A G$ such that $j^*\sigma_G^*
(\eta)=dx/x$, there exists a normal invariant differential
$\eta_x$ of $G_x$ such that $\eta=p_{G_x}^*(\eta_x)$.
\end{itemize}
\end{lemma}
\proof The first statement is immediate because
\begin{eqnarray*}
j^*f_x^*(\eta_x)=x^*\iota^*(\eta_x)=x^*(dz/z)=dx/x.
\end{eqnarray*}

For the second statement, let $\{S_h\}_h$ be an affine open covering
of $S$. For any $h$, let $\eta_{x,h}$ be a normal invariant
differential  of $G_x$ over $S_h$. The difference
$\omega_h=\eta-p_{G_x}^*\eta_{x,h}$ is the pull-back of an invariant
differential of $G$; moreover, as $p_{G_x}\circ\sigma_G\circ j=
f_x\circ j$, it holds
\begin{eqnarray*}
j^*\sigma_G^*p_{G_x}^*(\eta_{x,h})=j^*f_x^*(\eta_{x,h})=
dx/x=j^*\sigma_G^*(\eta).
\end{eqnarray*}
 Hence $\omega_h$ is indeed the pull-back of a suitable invariant
differential $\omega_{A,h}$ of $A_{S_h}$. Define $\tilde
\eta_{x,h}:=\eta_{x,h}+g^*\omega_{A,h}$. It satisfies
$p_{G_x}^*(\tilde \eta_{x,h})=\eta$, at least over $S_h$. Hence we
proved the assertion locally. To show that $\tilde \eta_{x,i}=\tilde
\eta_{x,h}$ on $S_i\cap S_h$ observe that $p_{G_x}^*\tilde
\eta_{x,i}= p_{G_x}^*\tilde \eta_{x,h}=\eta$ on $S_i\cap S_h$ and
$p_{G_x}^*\colon {\boldsymbol{\omega}}_{G_x}\to 
{\boldsymbol{\omega}}_{G_{x}\times_A
G} $ is injective. Hence the differentials $\tilde \eta_{x,i}$
provide a normal invariant differential $\eta_x$ of $G_x$ such that
$\eta=p_{G_x}^*\eta_x$. \qed

\proof(Proposition~\ref{pro.extnatex}) 
 By definition of $\delta$, it is 
$\alpha\circ \delta=0$. Moreover,
if $x$ is a character in ${\cal H}^\nabla(T)\cap {\cal H}(G)$ the extension
$G_x$ is isomorphic to the trivial one and the pull-back of $\eta$ to $A$ is
zero because it becomes zero on $G$ and 
${\boldsymbol{\omega}}_A\to {\boldsymbol{\omega}}_G$ has trivial kernel.
Let $(E,\nabla)$ be a $\natural$\nobd-extension of $A$ by the multiplicative
group. Suppose that its image via $\alpha$ is trivial. Hence we may
think $E$ as the push-out $G_x$ of $G$ with respect to a character
$x\colon T\to \mgr S$ and $\nabla$ as the connection associated to a
normal invariant differential $\eta_\nabla$ on $G_x$. It holds ${\bm
\omega}_{G_x}={\bm \omega}_{\mgr S}\times_{{\bm \omega}_T}{\bm
\omega}_{G}$. The projection of $\eta_\nabla$ on ${\bm \omega}_{\mgr
S}$ is $dz/z$ and the projection on ${\bm \omega}_{G}$ is  $du/u$
for a suitable homomorphism $u\colon G\to \mgr S$ (because the
pull-back of $(E,\nabla)$ to $G$ is isomorphic to the trivial
$\natural$-extension). Moreover the image of $dz/z$ in ${\bm
\omega}_T$ is $dx/x$ and it must coincide with the image of $du/u$
in  ${\bm \omega}_T$.   As the character $x/u_{|T}$ provides an extension
isomorphic to $G_x$ we may assume that $dx/x=0$. Hence $(E,\nabla)$ lies in
the image of $\bar\delta$.  
To show that $\beta\circ \alpha=0$, let $(E,\nabla)$ be a
$\natural$\nobd-extension of $A$ by the multiplicative group and
denote by $\eta_\nabla$ the normal invariant differential of $E$
associated to $\nabla$. Let $(E_G,\nabla_G)$ be the  pull-back of
$(E,\nabla)$ to $G$. Recall  that we have an exact sequence
\begin{eqnarray*} \xymatrix{0\ar[r] & {\boldsymbol{\omega}}_E\ar[r]^(.5){p_E^*} &
{\boldsymbol{\omega}}_{E_G}\ar[r] & {\boldsymbol{\omega}}_T \ar[r] & 0 }
\end{eqnarray*}
obtained as push-out of
\begin{eqnarray}\label{seq.invdiffAGT}
\xymatrix{0\ar[r] & {\boldsymbol{\omega}}_A\ar[r] & 
{\boldsymbol{\omega}}_{G}\ar[r] &
{\boldsymbol{\omega}}_T \ar[r] & 0. }
\end{eqnarray}
Let now $(E_T,\nabla_T)$ be the pull-back of $(E_G,\nabla_G)$ to
$T$. Clearly $E_T$ is isomorphic to the trivial extension and the
image of $p_E^*(\eta_\nabla)$ in  ${\boldsymbol{\omega}}_T$ is $0$; hence
$(E_T,\nabla_T)$ is isomorphic to the trivial
$\natural$\nobd-extension.

 We show now that the monomorphism $\coker(\bar\delta)\to
  \ker(\beta)$ is indeed an isomorphism.   
Suppose given a $\natural$\nobd-extension $(E_G,\nabla_G)$ of $G$ by
the multiplicative group over an $S$-scheme $S'$ such that its image via
$\beta$ is trivial.  We may assume that $E_G$ is the pull-back of an extension
$E$ of $A$ by $\mgr {S'}$ and that $S'$ is affine. Denote by $\eta_G$ the
normal invariant differential of $E_G$ associated to $\nabla_G$ and by 
$p_E\colon E_G\to E$ the projection homomorphism. The pull-back of $\eta_G$ to $T$ is an
invariant differential of type $dx/x$ for $x$ a character of $T$. As $S'$ is
affine, $E$ admits a $\natural$\nobd-structure $\nabla_E$ (associated to a
normal invariant differential $\eta_E$), so that $\omega=\eta_G -
p_E^*\eta_E$ is an invariant differential of $G$ and the restriction
of $\omega$ to $T$ is $dx/x$. The second statement of
Lemma~\ref{lem.normaldiff} asserts that the trivial extension with
the connection induced by $\omega$ is isomorphic to the pull-back of
$(G_x,\nabla')$ for a suitable connection $\nabla'$, hence its
isomorphism class lies in the image of $\alpha$. In particular
\begin{eqnarray*}
[(E_G,\nabla_G)]=\alpha[(E,\nabla_E)+(G_x,\nabla')]
\end{eqnarray*}
 and we get the result.

The exactness on the right can be deduced from
(\ref{seq.invdiffAGT}) and (\ref{seq.extnatM}) for $M=T, G$. 
\qed

In a similar way one gets the more general statement:
\begin{proposition}\label{pro.extnatexM} 
Let $M=[u\colon X\to G]$ be an $S$\nobd-$1$\nobd-motive. 
The following  sequence
\begin{eqnarray*}
\xymatrix
{
\uextnat {[X\to A]}{\mgr S}\ar[r]^(.5)\alpha
& \uextnat M{\mgr S}\ar[r]^(.5)\beta& \uextnat T{\mgr S}\ar[r]&0 }
\end{eqnarray*} is exact,  where $\alpha,\beta$ are the usual
pull-back morphisms and $\ker(\alpha)=
\displaystyle\frac{ {\cal H}^\nabla(T) }{ {\cal H}^\nabla(T) \cap {\cal H}(M) }$. 

If $S$ satisfies the hypothesis $(*)$,
$\alpha$ is a monomorphism. More generally, if $S$ is flat over $\ZZ$, the restriction of the above sequence to the site $({\bf Fl}/S)_{\rm fl}$ is also exact on the left.  
\end{proposition}
 
We will see in Corollary~\ref{cor.intersection} that ${\cal H}^\nabla(T) \cap {\cal H}(M)= {\cal H}^\nabla(M)$.

\subsection{Universal extensions and $\natural$-extensions.}

Recall that the universal extension $A^\natural$ of
an abelian variety $A$ represents  the sheaf $\uextnat{A'}{\mgr S}$
(cf. \cite{MM}). This does not extend to $1$\nobd-motives
in general.

\begin{proposition}
Let $M=[u\colon X\to G]$ be an $S$\nobd-$1$\nobd-motive and
$M^\natural=[u^\natural\colon X\to {\mathbf{G}}^\natural]$ its universal
extension.  
  There is a canonical epimorphism 
\[\psi_M\colon {\mathbf{G}}^\natural \longrightarrow   {\boldsymbol{\omega}}_{T'}
\times_{\uextnat{T'}{\mgr S}}\uextnat{M'}{\mgr S} \]
whose kernel is  $\ker(\alpha')=
\displaystyle\frac{ {\cal H}^\nabla(T') }{ {\cal H}^\nabla(T') \cap {\cal
    H}(M') }$. 
\end{proposition}

 \proof
 By the universal property of the push-out
we get from (\ref{dia.vectorial}) and  (\ref{dia.alphagamma}), for $G$ 
in place of $G'$, an epimorphism $\varphi_M$
 making the following diagram to commute
\begin{eqnarray}\label{dia.varphi}
\xymatrix {
 0\ar[r]&{\boldsymbol{\omega}}_{G'}\ar[r]^(.5)\iota \ar@{=}[d]&{\mathbf{G}}^\natural\ar[r]^(.5)\rho
 \ar@{->>}[d]^{\varphi_M} &
G\ar[r]\ar@{->>}[d]^\gamma &0\\
&{\boldsymbol{\omega}}_{G'} \ar[r]^(.3){\bar\iota}& \uextnat{M'}{\mgr S} \ar[r]
& \uext{M'}{\mgr S}\ar[r] &0.
}
\end{eqnarray}
We show that  $\varphi_M$ fits also in the following diagram
\begin{eqnarray}\label{dia.varphi2}
\xymatrix @C=20pt {
 0\ar[r]& G^\natural\ar[r]\ar@{=}[d] &{\mathbf{G}}^\natural\ar[r]^(.5)\tau
 \ar@{->>}[d]^{\varphi_M} & X\otimes \agr S={\boldsymbol{\omega}}_{T'}
\ar[r]\ar@{->>}[d]^j&0\\
 &  \uextnat{[X'\to A']}{\mgr S}\ar[r]^(.6){\alpha'} &\uextnat{M'}{\mgr S}
\ar[r]^(.5){\beta'} & \uextnat{T'}{\mgr S} \ar[r]&0 }
\end{eqnarray}
where the upper sequence in the one in (\ref{dia.vectorial}), the lower one
comes from Proposition~\ref{pro.extnatexM}, $j$ is the map in (\ref{seq.extnatM}).
To prove that $j\circ \tau=\beta\circ \varphi_M$ one proceeds as
follows:
 We may work locally and suppose that the vertical sequences
in (\ref{dia.vectorial}) are split.  Let $\delta$ be a section of
$\bar \tau\colon {\boldsymbol{\omega}}_{G'}\to {\boldsymbol{\omega}}_{T'}$. 
Any point
in ${\mathbf{G}}^\natural$ may be written as the sum $g+\iota(\delta(\omega))$
with $g$ a point of $G^\natural$ and $\omega$ a point of ${\bm
\omega}_{T'}$. Now,
\begin{eqnarray*}
\beta'(\varphi_M(g+\iota(\delta(\omega)) )
               &=&\beta'( \alpha'(g))+ \beta'(\bar\iota(\delta(\omega)))=
j(\bar \tau(\delta(\omega)))= j(\omega)~,\\
j(\tau(g+\iota(\delta(\omega))))&=& j(\tau(\iota(\delta(\omega))))=
                                   j(\bar \tau(\delta(\omega)))=j(\omega),
\end{eqnarray*}
because $\beta'\circ\alpha'=0$ by Lemma~\ref{pro.extnatexM} and 
$\beta'\bar\iota=j\bar\tau$.   
Diagram (\ref{dia.varphi2}) assures the existence of $\psi_M$
  whose kernel is isomorphic to the kernel of $\alpha'$. 
\qed

\begin{corollary}\label{cor.intersection}
With notations as above, it holds ${\cal H}^\nabla(T') \cap {\cal H}(M')={\cal
  H}^\nabla(M')$.
\end{corollary}
\proof Comparing all the previous constructions we get a cross of exact sequences
\begin{eqnarray*}
\xymatrix { &\ker{\psi_M}={\cal H}^\nabla(T')/{\cal H}^\nabla(T') \cap
  {\cal H}(M')\ar[d]\ar[rd]& \\
{\cal H}(M')/{\cal H}^\nabla(M')=\ker{\bar\iota}\ar[r]\ar[dr]& \ker{\varphi_M}\ar[r]\ar[d]& \ker\gamma={\cal
  H}(T')/{\cal H}(M')\\
 &\ker{j}={\cal H}(T')/{\cal H}^\nabla(T')& }
\end{eqnarray*}
where the upper diagonal arrow is a monomorphism by construction. Hence also
the lower diagonal arrow is a monomorphism and this happens if and only if 
 ${\cal H}^\nabla(T') \cap {\cal H}(M')={\cal H}^\nabla(M')$.
\qed

\begin{corollary}
\label{cor.EE}
Let $M$ be an $S$\nobd-$1$\nobd-motive with $S$ that satisfies hypothesis $(*)$.
The group scheme ${\mathbf{G}}^\natural$ in (\ref{dia.vectorial})  
represents the fibre product
\begin{eqnarray*}{\boldsymbol{\omega}}_{T'}
\times_{\uextnat{T'}{\mgr S}}\uextnat{M'}{\mgr S}.
\end{eqnarray*}
\end{corollary}

\begin{proposition}
Let $M$ be a $1$-motive over $S$ and $\varphi_M$ the epimorphism defined 
in (\ref{dia.varphi}). Once fixed a universal extension 
$[u_X^\natural\colon X\to {\boldsymbol{\omega}}_{T'}]$ of $[X\to
0]$ and ${\mathbf{G}}^\natural$ as in (\ref{dia.vectorial}), there
exists a canonical universal extension 
$M^\natural=[u^\natural\colon X\to {\mathbf{G}}^\natural]$ of $M$  
such that $\tau\circ u^\natural=u^\natural_X$ and the sequence 
\begin{eqnarray*}
\xymatrix{&  X \ar[r]^{u^\natural} & 
{\mathbf{G}}^\natural \ar[r]^(0.3){\varphi_M}& \uextnat{M'}{\mgr S}\ar[r] &0
}\end{eqnarray*} 
is exact. In particular, the kernel of $\varphi_M$ is isomorphic to ${\cal
  H}(T')/{\cal H}^\nabla(M')$. 
\end{proposition}
 
\proof
 
Uniqueness. Suppose  $u^\natural_1, u^\natural_2$ are universal extensions
such that $\varphi_M\circ u_i^\natural=0$ and 
$\tau\circ u_i^\natural=u_X^\natural$. 
Clearly $u^\natural_1- u^\natural_2$ factors through 
${\boldsymbol{\omega}}_{G'}$; as $\varphi_M\circ(u^\natural_1-
u^\natural_2)=0$, the morphism $u^\natural_1- u^\natural_2$ factors through
the subsheaf  $\ker \iota={\cal H}(M')/{\cal H}^\nabla(M')$ of 
${\boldsymbol{\omega}}_{G'}$. 
Furthermore $\tau\circ (u^\natural_1- u^\natural_2)\colon X\to
{\boldsymbol{\omega}}_{T'}$ is the zero map. 
It follows from Corollary~\ref{cor.intersection}, 
that the composition ${\cal H}(M')/{\cal H}^\nabla(M')\to {\cal H}(T')/{\cal
  H}^\nabla(T')\to {\boldsymbol{\omega}}_{T'}$ 
is a monomorphism; then $u^\natural_1=u^\natural_2$.

The uniqueness result assures that we can construct $u^\natural$ \'etale
locally. We proceed as in \cite{ABV}, 2.3, assuming that $X={\cal
  H}(T')=\oplus_i \ZZ e_i$. Let $\delta$ be a section of $\tau\colon
{\boldsymbol{\omega}}_{G'}\to {\boldsymbol{\omega}}_{T'}$ as in the proof of
Proposition~\ref{pro.repE} so that we identify ${\mathbf{G}}^\natural $ with
$G^\natural\oplus {\boldsymbol{\omega}}_{T'}$. If $\tilde u\colon X\to
G^\natural$ is a lifting of $u$, $u^\natural\colon X\to {\mathbf G}^\natural$
can then be defined via $u^\natural(e_i)=\tilde u(e_i)+\delta u_X^\natural(e_i)$.
 
  Recall that $u(e_i)$ is the $\mgr S$-extension $[X\to G_{e_i}]$ of $M'$ 
obtained as the push-out of $M'$ with respect to the character $e_i$. 
Let $f_i\colon G\to G_{e_i}$ be the induced map.
The section $u(e_i)$ lifts to a section $\tilde u(e_i) $ of
$G^\natural=\uextnat{M_A}{\mgr S}$ as soon as we fix an invariant
differential $\eta_i$ of $G_{e_i}$; locally this is always possible.  Then 
$\varphi_M\circ u^\natural (e_i)$ corresponds to the trivial $\mgr
S$-extension of $M$ together with the invariant differential
$f_i^*\eta_i+\delta u_X^\natural(e_i)$. Applying
Lemma~\ref{lem.normaldiff} it is immediate to check that  
$\eta_i$ can be chosen so that $\varphi_M\circ u^\natural(e_i)=0$.   
\qed

The previous proposition \emph{does not} imply that for any
universal extension $v\colon X\to {\mathbf{G}}^\natural$ of $M$ it holds $\varphi_M\circ
v=0$. Indeed, we have seen in Remark~\ref{rem.universal} that $v+f$
is also a universal extension for any homomorphism $f\colon X\to
{\boldsymbol{\omega}}_{A'}(\to {\boldsymbol{\omega}}_{G'}\to {\mathbf G}^\natural)$ and clearly $\varphi_M\circ f$ is not trivial in
general.

\begin{corollary}
Let $M$ be an $S$\nobd-$1$\nobd-motive  with $S$  that satisfies hypothesis $(*)$.
  Consider the homomorphism
 \begin{eqnarray*}
v\colon X\to {\mathbf{G}}^\natural=
{\boldsymbol{\omega}}_{T'} \times_{\uextnat{T'}{\mgr S}}\uextnat{M'}{\mgr S}
\end{eqnarray*}
whose projection to the first (resp.\  second) factor is
$x\mapsto  dx/x$ (resp.\ the $0$ map). It is
a universal extension of $M$. Moreover, there is an
exact sequence
\begin{eqnarray*}
\xymatrix{0\ar[r]& X\ar[r]^(.5){ v}&{\mathbf{G}}^\natural\ar[r]^(.3){\varphi_M}&
\uextnat{M'}{\mgr S} \ar[r]&0.}
\end{eqnarray*}
\end{corollary}

Observe that the right hand square in (\ref{dia.varphi}) is cartesian as soon
as $\uhomo{M'}{\mgr S}=0$, for example if $T'=0$. Under this hypothesis, we
could use the homomorphism $\varphi_M$ in (\ref{dia.varphi}) to prove the result in 
Proposition~\ref{pro.natpoincare}, i.e.\  the existence of a
$\natural$\nobd-structure on ${\cal P}^\natural$. In the general case however,
the homomorphism $\varphi_M$  looses information because, if we know the class
in $\extnat{M_{S'}}{\mgr {S'}}$ of a $\natural$\nobd-extension $({\cal
  P}_g,\nabla)$ with ${\cal P}_g$ the fibre of the Poincar\'e biextension of
$(M,M')$ at $g\in G'(S')$, we can determine $\nabla$ only up to an invariant 
differential of the type $du/u$ for $u$ a homomorphism of $M'_{S'}\to \mgr {S'}$.

{\bf Acknowledgements:}
 We thank F. Andreatta for a useful discussion on
Deligne's pairing. We are deeply indebted to V.~Cristante for helpful
comments on the first draft of the paper.  We thank the referee for pointing
out a problem in the use of a flatness hypothesis in the last section.

\end{document}